\documentclass[11pt]{article}
\usepackage{amsmath}
\usepackage{amssymb}
\usepackage{amsthm}
\usepackage{amsfonts}
\usepackage{amscd}
\usepackage[mathscr]{eucal}
\newcommand{\Z}{{\mathbb Z}}
\newcommand{\Q}{{\mathbb Q}}
\newcommand{\F}{{\mathbb F}}
\newcommand{\C}{{\mathbb C}}

\newcommand{\g}{{\mathfrak g}}
\newcommand{\h}{{\mathfrak h}}
\begin{document}
\parindent 25pt
\baselineskip 20pt
\textwidth 14.6cm
\textheight 22cm
\pagestyle{plain}

\title{ {
$L-$-Series and Their 2-adic and 3-adic Valuations at s=1 Attached
to CM   Elliptic Curves
 }
}
\author{\mbox{}\vspace{0.6cm}
{{\LARGE {QIU}} {\LARGE D}ERONG  \hskip 0.6cm
AND \hskip 0.6cm
 {\LARGE ZHANG} {\LARGE X}IANKE   } }
\date{}
\maketitle
\parindent 24pt
\baselineskip 18pt
\parskip 0pt

\par \vskip 1cm
\begin{center}
{\LARGE\bf   Abstract}
\end{center}
\vskip 30pt

$L-$series of  two classical families  of
elliptic curves  with complex multiplication are studied,
formulae for their special values at $s=1,\  $
bound of the values, and criterion of reaching the bound
are given. Let  $ E_{D}:\  y^{2}=x^{3}-D x \  $
be elliptic curves over the Gaussian field
$K=\Q(\sqrt{-1}), \ $
with $\ D=\pi _{1} \cdots \pi _{n}\  $  or
$\ D=\pi _{1} ^{2 }\cdots \pi _{r} ^{2} \pi _{r+1}
\cdots \pi _{n}$ ( where $\pi _{1},\ \cdots ,\ \pi _{n}$
are distinct primes in $K$. A formula for special values
of Hecke $L-$series  of such curves expressed by
Weierstrass $\wp-$function are given; a lower bound of
2-adic valuations of these values of
Hecke $L-$series  as well as a criterion for
reaching these bounds are obtained;
moreover the first part of the conjecture of Birch and
Swinnerton-Dyer is  hence
verified to hold for some of these curves.
Let $ E_{D^{2}}:\ \ y^{2}=x^{3}-2^{4}3^{3}D^{2} $
\quad and \quad
$ E_{D^{3}}:\ \ y^{2}= x^{3}+2^{4}D^{3}  \quad $
be elliptic curves over the field
$\, \Q(\sqrt{-3})\, $   \quad
(with $\ D=\pi _{1} \cdots \pi _{n}, \ $\,
where $\pi _{1} ,
\ \cdots ,\ \pi _{n}$ are distinct primes of
$\Q(\sqrt{-3})$), similar results as above
for $3-adic$ valuation are obtained.
These results develop some results
for more special case and for $2-adic$ valuation.
\par \vskip 1 cm

\begin{center}
{\Large {\bf   I. Introduction and Statement of Main Results}}\\
\end{center}

Consider the two classical families  of elliptic curves:\
$$E_1:\ \quad  y^2=x^3 - D x ; $$
$$E_2:\ \quad   y^2=x^3 + D^\prime ,  $$
here $D ,\ D^\prime  \in \mathbb{Z} $ are rational integers
(but in the  following we will generalize to the cases
that $D ,\ D^\prime  $ are certain quadratic
algebraic numbers).
These elliptic curves
have been studied broadly for a long time,
having relations with
many problems of number theory. For example, the curve
$ E_1 $ correlates intimately  to the problem of
congruent numbers when  $D $
is a square in $\mathbb{Z}$(see [Tun] ).
These two families  of elliptic curves $E$ have
complex multiplication by $\sqrt{-1}$ and $\sqrt{-3}$
respectively, their complex $L-$series (or $L-$function)
$L(E,\ s)$ could be
identified with the $L-$series attached to  certain
Hecke characters (i.e. Gr$\ddot{o}$ssencharacter) of
the fields ${\Q}$ $ (\sqrt{-1})$ and $ {\Q }$ $(\sqrt{-3})$
respectively.
The ``conjecture of Birch and  Swinnerton-Dyer"
(or `` B-SD conjecture", for brevity) asserts that
the value  of the  $L-$series $L(E, \ s)$ at $s=1$ of an
elliptic
curve $E,\ $  $L(E, \ 1),\ $   is very important for the
arithmetic study of
the elliptic curve. A considerable of (numerical) evidences
for the B-SD conjecture have been held up since
it was published,  most of them were
from the above two families of elliptic curves
as could be seen in the original paper of
Birch and Swinnerton-Dyer [B-SD] , and papers of
Razar [Razar ??]  and Stephens [Ste] .
In particular, if $D $ is a perfect square of
integer, then the B-SD conjecture predicts that
$L({E_1},\ 1)$ could be divided by  a power of $2$
(up to a multiple of an appropriate period of $E_1$) which
depends on the number of distinct prime
factors of $D $. For certain kinds of $D^\prime ,\ $
the B-SD conjecture has similar prediction for
the curve  $ E_2 $
(see e.g. [Razar] , [Tun] , [Ste] ).
\par \vskip 0.3cm

In 1997,  C. Zhao studied the problem of
divisibility of $L({E_1},\ 1)$  by  powers of $2$
under the new assumption $D  = D_0^2$ with
$D_0$ an Gaussian integer in the Gaussian field
$\Q (\sqrt{-1})$. In fact, he studied
the $2-adic$ valuation of the value  at $s=1$ of
the $L-$series of the elliptic
curve $E_1:\ \quad  y^2=x^3 - D_0^2x ,
\ D_0\in \Z [\sqrt{-1}] $
(Actually, the value of the $L-$series should be divided
by an appropriate period $\Omega $ of the elliptic curve first;
this normally will not be mentioned again  in the following
as a default fact);
he gave the  rigorous lower bound for the
$2-adic$ valuation as well as a criterion of
reaching this bound,  and hence obtained
nice results about congruent numbers and
showed the first part of the B-SD conjecture
is true for some elliptic curves $E_{D_0}^2$ over
Gaussian field.
\par \vskip 0.3cm

In 1968, N. Stephens studied a case of the  elliptic curves
$E_2,\ $  i.e.  $E:\  \ y^2= x^3 - 2^43^3D_1^2\ $ with
$D_1\in \Z$([Ste] ). He proved that
if $D_1>2$ is a cube-free rational integer, then
$\ \psi (D_1)= 3^{1/2}D_1^{1/3}L(E, \ 1)/\Omega \ $
is a rational integer (where $\Omega $ is the period
mentioned above, a constant expressed by Weierstrass
$\wp -$functions), and
\par \vskip 0.2cm
$3$ divides $ \psi (D_1)= 3^{1/2}D_1^{1/3}L(E, \ 1)/\Omega $
\qquad (when $9|D_1) .$  \par \vskip 0.2cm

In the present paper, we will study the two classical
families  of elliptic curves $E_1$ and $E_2 $
further on the base fields $K_1={\Q} (\sqrt{-1})$
and $K_1={\Q} (\sqrt{-3})$ respectively,
giving formulae of values of their $L-$series at $s=1,\ $
lower bound for their
$2-adic $ and $3-adic $ valuations ,
criteria for reaching the bounds,
and  verify  B-SD conjecture in some cases.

\par \vskip 0.3cm

Over the Gaussian field $\Q(\sqrt{-1}),\ $  we will first
study the elliptic curves
$ E_{D}:\ y^{2}=x^{3}-Dx \quad  $
with $\ D=\pi _{1} \cdots \pi _{n}\ $ and
$\ D=\pi _{1} ^{2 }\cdots \pi _{r} ^{2} \pi _{r+1}
\cdots \pi _{n}\ $
(where $ \ \pi _{1},\ \cdots ,\ \pi _{n}\ $
are distinct Gaussian prime integers in $\Z(\sqrt{-1})\ $
(when $r=n,\ $ the second case turns to be the case
studied by Zhao [Zhao] ). We will give a formula for
the special values at $s=1$ of the
Hecke  $L-$Series of $E_D$ ( expressed via Weierstrass
$ \wp -$function), lower bounds for the
$2-adic$ valuation of the values, criterion of reaching
the bounds, and show that the B-SD conjecture about the
relation of the rank of the Mordell-Weil group
and the analytic rank of the $L-$series is true for some
elliptic curves $E_D,\ $  by using our criterion and results
of Coates and Wiles .
\par \vskip 0.4cm

Then over the quadratic field
$ \Q(\sqrt{-3}),\ $  we consider the two kinds of
elliptic curves
$ E_{D^{2}}:\ \quad y^{2}=x^{3}-2^{4}3^{3}D^{2},  $
and $ E_{D^{3}}:\ \quad y^{2}=
x^{3}+2^{4}D^{3} \quad  $
with $\quad D=\pi _{1} \cdots \pi _{n}\, $ where
$\ \pi _{1} ,\ \cdots ,\ \pi _{n}\ $ are distinct prime
integers in $K_2= \Q(\sqrt{-3})$. Similar results
as above (but for $3-adic$ valuation) will be given.
These results develop the results about the estimation of
$2-adic$ valuation for special value of $L-$series of $E_D$
with $D=D_0^2 $ square in Gaussian field in [Zhao] .
\par \vskip 0.6cm

(A.1) Now let   $K=\Q(\sqrt{-1})$ be the Gaussian field, $O_K=\Z(\sqrt{-1})$  the ring of its integers
(Gaussian integers), and put $I=\sqrt {-1}$.  Consider the
elliptic curve
$$E_{D}:\  \ y^{2}=x^{3}-Dx \hskip 1cm
D=\pi _{1} \cdots \pi _{n},\hskip 2cm  $$
where
$\pi _{1},\ \cdots ,\  \pi _{n}$
are distinct prime (Gaussian ) integers in
$O_K=\Z(\sqrt{-1})$  and
$\pi _{k}\equiv 1(mod \  4)$
($k=1,\ \cdots,\ n \ ) .\ $
So $D\equiv 1 (mod \  4) \ .$ \par
\vskip 0.3cm

Denote the set $S=\{\pi _{1},\ \cdots, \ \pi _{n}\}$.
For any subset $ T$ of the set $\{1,\ \cdots ,\ n\},\ $
define
$$ D_{T}=\prod _{k \in T} \pi _{k},\qquad
\widehat{D}_{T}=\prod _{k \not \in T} \pi _{k}
=D/D_{T} \ ,\hskip 2cm  $$
and put $D_{\emptyset}=1$ when $T=\emptyset $ (emptyset).
Let
$$  L_{S}(\overline{\psi } _{D_{T}}, \ s)$$
denote the Hecke $ L-$series of $\psi _{D_{T}}$
(omitting all the Euler factors corresponding to
primes in $S$ ), where $\psi _{D_{T}}$ is the
Hecke character (Gr$\ddot{o}$ssencharacter)
of the field $K$ corresponding to the elliptic curve
$E_{D_{T}}:\ \ y^{2}=x^{3}-D_{T}x$ . For the
special value $L_{S}(\overline{\psi }_{D_{T}},\quad 1)$
of the above $L-$series at $s=1,\ $  we have the following
formula  expressed as a finite sum of
Weierstrass $ \wp -$function $\wp (z)$.
\par \vskip 0.5cm

{\bf Theorem 1.1}  \quad  For any factor $D_T$ of
$D=\pi _{1} \cdots \pi _{n}\in \Q (\sqrt {-1}) $ as above ,
let $\psi _{D_{T}}$ be the Hecke character of the Gaussian
field $\Q (\sqrt {-1}$  corresponding to the elliptic curve
$E_{D_{T}}:\ \ y^{2}=x^{3}-D_{T}x$.
Then we have
$$ \frac{D}{\omega }
\overline{\left(\frac{\theta}{D_{T}}\right)_{4}}
L_{S}\left(\overline{\psi }_{D_{T}},\ 1\right)=\frac{I}{2} \sum _{c \in \mathcal{C}}
\left(\frac{c}{D_{T}}\right)_{4} \frac{1}{\wp \left(\frac{c\omega }{D}\right)-I}+
\frac{1}{4}\sum _{c \in \mathcal{C}} \left(\frac{c}{D_{T}}\right)_{4}
\qquad (1.1) $$
where $\theta =2+2I,\ \left(\frac{\alpha }{\beta }\right)_{4} $
denote the (generalized) quartic residue symbol, $\mathcal{C}$ is any complete reduced
remainder-system of $O_{K}$ modulo $D,\ $
$L_{\omega }=\omega O_{K}$ is the period lattice of the
elliptic curve $E_{1}:\ y^{2}=x^{3}-x,\ $
$$\omega = \int _{1} ^{\infty } \frac{dx}{\sqrt{x^{3}-x}}
=2.6220575 \cdots , $$
$ \wp (z)$ is the Weierstrass $\wp-$function associated to
the lattidce   $L_{\omega }$.
\par \vskip 0.5cm

Let $\overline{\Q _{2}}$ be the algebraic closure of
the $2-adic$ (complete) field $\Q _{2},\ $
$\ v=v_2 $ is the normalized $2-adic $ exponential
valuation of
$\overline{\Q _{2}}\ $ (i.e. $v_2(2)=1  ).\ $
Fix an isomorphic embedding $\overline{\Q} \hookrightarrow
\overline{\Q _{2}},\ $ where $\overline{\Q}$ is the algebraic
closue of the rational field $\Q$. For any algebraic
number $\alpha ,\ $ let $v_2(\alpha )$ denote the
$2-adic $ valuation of $\alpha $.
For $D=\pi _{1} \cdots \pi _{n}$ as above, put
$$ S^{*}(D)=\frac{I}{2} \sum _{c \in \mathcal{C}}
\frac{1}{\wp (\frac{c\ \omega}{D})-I} \sum _{T}
\left(\frac{c}{D_{T}} \right) _{4}.  \qquad (1.2)  $$
For any Gaussian integers $\alpha ,\ $ $\beta \ $ which
are relatively prime, put
$\quad \left( \frac{\alpha }{\beta } \right) _{4} ^{2}=
\left( \frac{\alpha }{\beta } \right) _{2} \ , $ and define
$\ \left[ \frac{\alpha }{\beta } \right]  _{2}=
(1-\left( \frac{\alpha }{\beta } \right) _{2})/2\
,\ $ then
$\quad \left[ \frac{\alpha \gamma }{\beta } \right]  _{2}=
\left[ \frac{\alpha }{\beta } \right]  _{2}+
\left[ \frac{\gamma }{\beta } \right]  _{2}$ \\
(regard $  [ - ]  _{2}$ as a $\F_{2}-$value function).
\par \vskip 0.3cm

For $D=\pi _{1} \cdots \pi _{n},\ $ we define a
$\F_{2}-$value function $\delta _k $ inductively
as the following :
First put
$$ \varepsilon _{n}(D)=\left\{
\begin{array}{l}  1,  \qquad \hbox{if}
\quad v_2(S^{*}(D))=(n-1)/{2}
\quad ; \\
0,  \qquad \hbox{if} \quad v_2(S^{*}(D))>
(n-1/{2}  \ ,
\end{array}
\right. $$
where $n=n(D) $ is the number of distinct prime factors of
$ D$. Then for Gaussian prime integer
$\pi \ $ with  $\pi \equiv 1 (mod \  4),\ $ define
$s_{1}$ as a  $\F _{2}-$valued function as follows:
$$ s_{1}(\pi )=\left\{
\begin{array}{l} 1 , \qquad \hbox{if} \quad
v_2(\pi -1)=2 \quad ;\\
0, \qquad \hbox{if} \quad v_2(\pi -1)> 2 \quad .
\end{array}
\right . $$
Finally define the $\F_{2}-$value function $\delta _{k}
\quad (k=1,2,\cdots )$ as follows:
$$ \delta _{1} (\pi )=s_{1}(\pi )+\varepsilon _{1}(\pi );  $$
and for $\quad D=\pi _{1} \cdots \pi _{n} \quad  (n\geq 2),\ $
define
$$ \delta _{n} (D)=\delta _{n} (\pi _{1},\cdots , \pi _{n})=
\varepsilon _{n}(D)+\sum _{\emptyset \neq T \subsetneqq
\{1,\cdots , n \} } \left(\prod _{k \not\in T}
\left[\frac{D_{T}}{\pi _{k}} \right]  _{2} \right)
\delta _{t}(D_{T}) \ , $$
where the sum `` $\sum$'' is taken over the nonempty subsets $T$
of $\{1,\cdots , n\}$, and  $t=\sharp T$ is
the cardinal of $T$.
\par \vskip 0.5cm

{\bf Theorem 1.2.}\quad  Let $D=\pi _{1} \cdots \pi _{2}, \ $
where $\ \pi _{k}\equiv 1(mod \ 4)\ $ are distinct  Gaussian
prime integers  $(k=1,\cdots , n)$. Then for the
$2-adic $ valuation of the values of the $L-$series
we have
$$  v_2 \left(L(\overline{\psi } _{D},1)/\omega
\right)\geq \frac{n-1}{2},\ $$
and the equality holds if and only if
$\  \delta _{n} (D) = 1 \  $ .
\par \vskip 0.5cm

{\bf Theorem 1.3.} \quad Let $D=\pm p_{1} \cdots  p_{m}
\equiv 1 (mod \ 4),\ $
where $\quad p_{k} \not\equiv 5 (mod\ 8)\ $
are distinct positive rational prime numbers
$(k=1,\cdots , m ). $ If $\ \delta _{n}(D)=1,\ $
then the first part of B-SD conjecture  is true for the
elliptic curve $E_{D}:\ y^{2}=x^{3}-Dx$, that is
$$rank(E_{D}(\Q))= ord_{s=1}( L(E_{D}/\Q , s))=0 .
\hskip 1.5cm  (1.3)$$
(where $n$ is the number of distinct prime factors of $D$).
\par \vskip 0.5cm

(A.2) Now consider the elliptic curves
$\quad E_{D}:\ y^{2}=x^{3}-Dx \ $
for $\ D =\pi _{1} ^{2} \cdots \pi _{r} ^{2} \pi _{r+1}
\cdots \pi _{n},\ $
where $\quad \pi _{k}\equiv 1 (mod\ 4)$ are distinct prime
Gaussian integers $(k=1,\ \cdots ,\ n).\ $
Similarly to tht above, let
$\quad S=\{ \pi _{1},\ \cdots ,\ \pi _{n}\}.\ $
Write any subset $T$ of $\{ 1,\ \cdots ,\ n\}$
as $T=T_{1}\cup T_{2},\ $  where
$T_{1}=T \cap \{1,\ \cdots ,\ r \} ,  \qquad
T_{2}=T\cap \{ r+1,\ \cdots ,\ n \}  .$
And define
$$ D_{T}=D_{T_{1}}D_{T_{2}}\qquad ,\quad
D_{T_{1}}=\prod _{k \ \in T_{1}}
\pi _{k} ^{2}\qquad ,\quad D_{T_{2}}=
\prod _{k \in T_{2}} \pi _{k} \ .$$
$$ \widehat{D}_{T}=\frac{D}{D_{T}}=
\widehat{D}_{T_{1}} \widehat{D} _{T_{2}}
 \qquad  D_{1}=\pi _{1} ^{2} \cdots \pi _{r} ^{2}  \ ,  \quad
D_{2}=\pi _{r+1} \cdots \pi _{n}  \ ,  \quad D=D_{1}D_{2}. $$
$$ \widehat{D}_{T_{1}}=\frac{D_{1}}{D_{T_{1}}}  \ ,  \quad
\widehat{D}_{T_{2}}=\frac{D_{2}}{D_{T_{2}}}  . $$
When $\quad T=\emptyset $(emptyset)( or $T_{1}=\emptyset ,\ $
or $T_{2}=\emptyset $), define $\quad D_{T}=1 \quad
$(or $D_{T_{1}}=1,\ $ or $D_{T_{2}}=1$ respectively).
(If $r=n$, then $D_{1}=D,\ $ $D_{2}=1$;
and if $r=0$, then $D_{1}=1,\ $ $D_{2}=D$). \\
For above  $\ D =\pi _{1} ^{2} \cdots \pi _{r} ^{2} \pi _{r+1}
\cdots \pi _{n},\ $  set $S$ and $T,\ $
denote
$$ \Delta  _{1}=\pi _{1} \cdots \pi _{r}  \qquad  \Delta  _{2}=
\pi _{r+1} \cdots \pi _{n}  \qquad
\Delta  =\pi _{1} \cdots \pi _{r}
\pi _{r+1} \cdots \pi _{n}=\Delta  _{1}\Delta  _{2}  \ ,  $$
$$ \Delta  _{T}=\prod _{k \in T_{1}} \pi _{k}\cdot \prod _{k \in T_{2}}
\pi _{k}=\Delta  _{T_{1}} \Delta  _{T_{2}}  \ ,  \quad \Delta  _{T_{1}}=
\prod _{k \in T_{1}} \pi _{k}  \ ,  \quad \Delta  _{T_{2}}=
\prod _{k \in T_{2}} \pi _{k}  . $$
$$ \widehat{\Delta  }_{T}=\frac{\Delta  }
{\Delta  _{T}}  \ ,  \quad
\widehat{\Delta  }_{T_{1}}=\frac{D_{1} }
{\Delta  _{T_{1}}}  \ ,  \quad
\widehat{\Delta  }_{T_{2}}=\frac{D_{2} }
{\Delta  _{T_{2}}}  . $$
Define $\Delta  _{\emptyset}=1$.
\par \vskip 0.3cm

Let $\quad L_{S}(\overline{\psi } _{D_{T} },\ s) \quad $
denote Hecke $  L-$series of $\psi _{D_{T} }$
(omitting all Euler factors corresponding to primes in
$S$), where $\psi _{D_{T} }$ is the Hecke character of the
Gaussian field $K=\Q (\sqrt{-1})$ corresponding to the
elliptic curve
$\quad E_{D_{T}}:\ y^{2}=x^{3}-D_{T}x \ $ .
By the definition we know
$$ L_{S}(\overline{\psi } _{D_{T} },\ s)=\left\{
\begin{array}{l}
L(\overline{\psi } _{D },\ s)  \qquad  \hbox{if } \quad
D_{T} =D \quad ;\\
L(\overline{\psi } _{D_{T} },\ s) \prod \limits _{\pi _{k}
|\widehat{D}_{T}} \left( 1-\left( \frac{D_{T}}{\pi _{k}}
\right) _{4}
\frac{\overline{\pi }_{k}}{(\pi _{k}
\overline{\pi } _{k}) ^{s}}
\right)  \qquad  \hbox{otherwise}.
\end{array}
\right. $$
\par \vskip 0.5cm

{\bf Theorem 1.4.} \quad  For any factor $D_T=D_{T_1}D_{}T_2$
of $\ D =\pi _{1} ^{2} \cdots \pi _{r} ^{2} \pi _{r+1}
\cdots \pi _{n}\in \Q (\sqrt {-1} $ as above ,
let $\psi _{D_{T}}$ be the Hecke character of the Gaussian
field  corresponding to the elliptic curve
$E_{D_{T}}:\ \ y^{2}=x^{3}-D_{T}x$.
Then we have
$$ \frac{\Delta  }{\omega } \overline{\left( \frac{\theta }
{D_{T}} \right) _{4}} L_{S}(\overline{\psi } _{D_{T} },\ 1)=
\frac{I}{2} \sum _{c \in \mathcal{C}}\left( \frac{c}
{D_{T}} \right) _{4} \frac{1}{\wp \left( \frac{c\ \omega }
{\Delta  }
\right)-I}+\frac{1}{4} \sum _{c \in \mathcal{C}}
\left( \frac{c}
{D_{T}} \right) _{4},  \qquad (1.)$$

where $\Delta=$ $??????????????????????????????????$
\par \vskip 0.5cm

$\mathcal{C}$ is a complete reduced
residue system of $O_{K}$ modulo $\Delta  $; $\theta =2+2I,\ $
$L_{\omega }=\omega O_{K}$, $\omega ,\ $  and $\wp (Z)$ are as in Theorem 1.1.
\par \vskip 0.5cm

{\bf Theorem 1.5.} \quad Let $D=\pi _{1} ^{2} \cdots
\pi _{r} ^{2} \pi _{r+1} \cdots\pi _{n},\ $
$n,\ r$ are positive integers $(1\leq r \leq n)$
(If $r=n$ then $D=\pi _{1} ^{2} \cdots \pi _{n} ^{2}$), $\pi _{k}\equiv 1(mod\ 4),\ $ are distinct prime Gaussian
integers ($k=1,\ \cdots ,\ n$). Then for the
$2-adic $ valuation of the values of the $L-$series we have
$$\quad v \left( L(\overline{\psi } _{D},\ 1)/\omega
\right) \geq \frac{n}{2}-1  .\hskip 1.5cm (1.5) $$
       \par \vskip 0.5cm

(B) Now we consider the elliptic curves
$ y^{2}=x^{3}-D^\prime $ over the number field
$K=\Q(\sqrt{-3})$
with complex multiplication by
$\sqrt {-3}$.
Let  $\tau =({-1+\sqrt{-3}})/{2}=exp({2\pi I/3})$ be
a primitive cubic root of unity, $O_{K}=\Z [\tau
] $ be the ring of integers of $K$.
We will study elliptic curve
$$E_{D^{2}}:\quad  y^{2}=x^{3}-2^{4}3^{3}D^{2}\ , $$
where $D=\pi _{1}\cdots \pi _{n},\ $
$\pi _{k}\equiv 1(mod\ 6)\ $ are distinct prime
elements of $O_K$ ( $k=1,\cdots , n 0.$
\par \vskip 0.3cm

Let $\ S=\{\pi _{1},\cdots , \pi _{n}\}.\ $ For any
subset $T$ of $\{1,\cdots  , n\}$,  define
$$ D_{T}=\prod _{k\in T} \pi _{k} \qquad ,\quad
\widehat{D} _{T}=
\prod _{k\not\in T}\pi _{k}={D}/{D_{T}} \  , $$
and put $D_{\emptyset}=1 \ .$  Let
$\psi _{D_{T} ^{2}}$ be the Hecke character of
of $K$ corresponding to the elliptic curve
$$E_{D_{T} ^{2}}:\ y^{2}=x^{3}-2^{4}3^{3}D_{T}
^{2}. $$
And let $L_{S}(\overline{\psi }
_{D_{T} ^{2}},\ s)$ denote the Hecke $L-$series
of $\psi _{D_{T} ^{2}}$ (omitting all
the Euler factors corresponding to primes in $S$).
Then
$\quad L_{S}(\overline{\psi } _{D_{T} ^{2}},\ 1)$
could be expressed by the
 Weierstrass $  \ \wp -$functions as in the following:
\par \vskip 0.5cm

{\bf Theorem 1.6.}\quad For any factor $D_T$ of
$D=\pi _{1} \cdots \pi _{n}\in \Q (\sqrt {-3}) $ as above ,
let $\psi _{D_{T}^2}$ be the Hecke character of the
field $\Q (\sqrt {-3})$ corresponding to the elliptic curve
$E_{D_{T} ^{2}}:\ y^{2}=x^{3}-2^{4}3^{3}D_{T}
^{2}. $ Then we have
$$ \frac{D}{\omega }\left(\frac{9}{D_{T}} \right) _{3} L_{S}(\overline
{\psi } _{D_{T} ^{2}},1)=\frac{1}{2\sqrt{3}} \sum _{c\in \mathcal{C}}
\left(\frac{c}{D_{T}} \right) _{3} \frac{1}{\wp \left(\frac{c\omega }{D} \right)
-1}+\frac{1}{3\sqrt{3}}\sum _{c\in \mathcal{C}}
\left(\frac{c}{D_{T}} \right) _{3} .\quad (1.6 ) $$
where $\wp (z)$ is the Weierstrass $\\wp -$function
associated to $L_{\omega }=\omega O_{K}$
the period lattice of the elliptic curve
$\quad E_{1}:\ y^{2}=x^{3}-\frac{1}{4}\quad $, $\mathcal{C}$ is a complete residue system of $O_{K}$
modulo $D$, $\omega =3.059908 \cdots \ $ is a constant, $\left(\frac{a}{b} \right) _{3} $ is the cubic residue
symbol.
\par \vskip 0.5cm

Now let $\overline{\Q}_{3}$ be the algebraic closure
of $\Q_{3}$, the  $3-adic $completion of $\Q$.
Let $v_{3}$ be the normalized ($3-adic $) exponential
valuation of
$\overline{\Q}_{3}$ , i.e.
$v_{3}(3)=1 \quad .$ Fix an (isomorphic) embedding
$\overline{\Q}\hookrightarrow \overline{\Q_{3}}\ .$
\par \vskip 0.5cm

{\bf Theorem 1.7.}
\quad Let $D=\pi _{1} \cdots \ \pi_{n},\ $
where $\quad \pi _{k}\equiv 1(mod\ 6),\ $
are distinct prime elements of
the number field $K=\Q(\sqrt{-3})$
($k=1,\ \cdots ,\ n \quad $). Then we have
$$ v_{3} \left(L(\overline{\psi } _{D^{2}},\ 1)/
\omega \right)\geq   \frac{n}{2}-1  .\hskip 3 cm (1.7) $$
\par \vskip 0.4cm

{\bf  Theorem 1.8.}
For $\quad D=\pi _{1}\cdots \pi _{n} \quad $
as in Theorem 1.7, put
$$ S^{*}(D)=\frac{1}{2\sqrt{3}} \ \sum _{c\ \in \mathcal{C}}
\frac{1}{\wp (\frac{c\ \omega }{D},\ L_{\omega })-1}
\sum _{T} 2^{n-t(T)}\ \left( \frac{c}{D_{T}} \right) _{3}.
\hskip 2cm (1.8) $$
Then we have
$$  v_{3}(S^{*}(D))\geq \frac{n-1}{2}  . $$
\par \vskip 0.6cm

\begin{center}
{\Large {\bf   II. 2-Valuations of $L-$series of Elliptic
Curves with CM by $\sqrt {-1}$ }}\\
\end{center}

We need the following results :
\par \vskip 0.4cm

{\bf  Proposition A}\quad Let $E$ be an elliptic curve defined over the imaginary
quadratic field $K$ with complex multiplication ring $O_{K}$( integers in $K$ ). Assume
its period lattice is
 $L=\Omega O_{K}$, $\Omega \in \C^{\times }$ a
 complex number, $\phi $ is the Hecke character of $K$
corresponding to $E$, $\g$ is an integral ideal of $K$
, $E_{\g}$ is the subgroup of $E$ consist of $\g-$divisible
points. Let $\mathbf{B}$ be a set of integral ideals of $K$
relatively prime to $\g$  and
$$ \{\sigma _{\flat } \ | \ \flat \in \mathbf{B} \}
=Gal(K(E_{\g})/K),\quad
(\hbox{if}\ \flat \neq \flat ',\ \hbox {then} \
\sigma _{\flat }\neq
\sigma _{\flat '} )$$
where
$$\quad \sigma _{\flat }=\left(\frac{K(E_{\g} )/K}{\flat}\right)$$
is $Artin$ symbol.
Put $\rho \in \Omega K^{\times }\subset \mathcal{\C}^{\times }$
, and $\rho \Omega ^{-1}
O_{K}=\g^{-1}\h$, $\h$ is an integral ideal of $K$
relatively priime to $\g$ . Then
$$ \frac{\phi ^{k}(\h)}{N(\h)^{k-s}} \cdot
\frac{ \overline{\rho } ^{k}}{|\rho | ^{2s}}
\cdot L_{\g}(\overline{\phi } ^{k},\quad s)=
\sum _{\flat \in \mathbf{B}} H_{k}
(\phi (\flat )\rho,\ 0,\ s,\ L) \  $$
($Re(s)>1+k/2 \ $),  where $k$ is an positive integer, $N$ denotes the norm map from $K$ to $\Q$ .
$$ L_{\g}(\overline{\phi } ^{k},\ s)=
\prod _{\wp \dag \g} (1-\overline{\phi }
^{k}(\wp )N(\wp ) ^{-s}) ^{-1}, \qquad (Re(s)> 1+k/2 )$$
$$ H_{k}(z,\ 0,\ s,\ L)= \sum{^{'}} \frac{
(\overline{z}+\overline{\alpha}) ^{k}}{|z+\alpha
| ^{2s}},\qquad (Re(s)>1+k/2) $$
where the sum  $\sum '$ is taken over
 $ \alpha \in L=\Omega O_{K}$ and $\alpha \neq -z $
  when $z \in L$
[Go-Sch]  .
\par \vskip 0.4cm

{\bf Lemma B. } \quad  Let elliptic curve $E$,
field $K$, Hecke character $\phi $, and
$\g$ are as in Proposition A. If the conductor
 $f_{\phi }$ of $\phi $  divides  $\g$, then
$K(E_{\g})$ is the ray class field of $K$
to the  cycle (or divisor, modulo ) $\g$
(see [Go-Sch] ).
\par \vskip 0.4cm

Now we consider Theorem 1.1 and let $K=\Q(\sqrt{-1})$ , $E_{D},\ $  $ D_{T}$,  and
$L_{S}(\overline{\psi } _{D_{T}}, \ s)$ etc
be as there. Then by definition
(see [B-SD] , [Ire-Ro] ) we have
  \par \vskip 0.4cm

{  \bf  Lemma 2.1.}
$$L_{S}(\overline{\psi } _{D_{T}}, \ s)=\left\{
\begin{array}{l} L(\overline{\psi } _{D_{T}}, \ s) ,
\hskip 2cm  \hbox{if} \quad
\prod\limits _{\pi _{k} \in S} \pi _{k}=D_{T} \quad ; \\
L(\overline{\psi } _{D_{T}}, \ s)\prod
\limits _{\pi _{k}|\widehat{D}_{T}}
(1-\left( \frac{D_{T}}{\pi _{k}} \right)_{4}\cdot
\frac{\overline{\pi _{k}}}
{(\pi _{k} \overline{\pi _{k}}) ^{s}}), \hskip 1cm
\hbox{otherwise . }
\end{array}
\right . $$
\vskip 0.5cm

{\bf  Proof of Theorem 1.1.}\quad
For the elliptic curve $E_{D_{T}}:\ y^{2}=x^{3}-D_{T}x$,
assume its period lattice is $L=\Omega O_{K},\ $ with
$\Omega =\alpha \omega,\ $ $\alpha \in \C^{\times }$
(Obviously $\Omega ={\omega }/{\sqrt[4] {D_{T}}}$).
From [Bir-Ste]  we know the conductor of $\psi _{D_{T}}$
is $(\theta D_{T}).\ $
Now, in Proposition A, let $k=1,\ $
$\rho =\Omega /(\theta D),\ $
$\g=(\theta D),\ $ $\h=O_{K}$, we have
$$ \frac{\overline{\rho }}{|\rho | ^{2s}}L_{\g}
(\overline{\psi }_{D_{T}}, \ s)=
\sum _{\flat \in \mathbf{B}} H_{1}(\psi _{D_{T}}(\flat )
\rho ,0, s, L), \quad (Re(s))> {3}/{2}).\hskip 1cm (2.1)  $$
Since the conductor of $\psi _{D_{T}}$ is $\theta D_{T},\ $
and $(\theta D_{T}) \ | \ (\theta
D)=\g,\ $ so by Lemma  B we know that the ray class field
of $K$ to the cycle $(\theta D)$ is
$K((E_{D_{T}})_{(\theta D)}),\ $ in particular we have
the following isomorphism via Artin map:
$$ (O_{K}/(\theta D)) ^{\times }/\mu _{4}
\cong Gal(K((E_{D_{T}}) _{(\theta D)})/K), $$
where $\mu _{4}$ is the group of quartic roots of unity,
and $\mu _{4}\cong (O_{K}/\theta ) ^{\times } \quad . $
So we could take the set
$$\mathbf{B}=\{(c\theta +D) \quad |
\quad c \in \mathcal{C} \},\hskip 3cm (2.2) $$
where $\mathcal{C}$ are fixed representations of
$(O_{K}/(D)) ^{\times }$, so we have
$$ \frac{\overline{\rho }}{|\rho | ^{2s}}L_{\g} (\overline{\psi } _{D_{T}}, \ s)=
\sum _{c \in \mathcal{C}} H_{1}(\psi _{D_{T}}(c\theta +D)\rho ,0, s, L),
\quad (Re(s)>{3}/
{2})  \hskip 1cm (2.3) $$
Note that the analytic extension of $H_{1}(z,o,1,L)$
could be obtained by Eisenstein $\ E^{*}-$
(see [Zhao]  or [We] ), i.e., $H_{1}(z,0,1,L)=E^{*}_{0,1}(z,L)=E_{1}^{*}(z,L) \ , $
So by (2.3)  we have
$$ \frac{\theta D}{\Omega }L_{(\theta D)}
(\overline{\psi }_{D_{T}} ,\quad 1)=
\sum _{c \in \mathcal{C}}E_{1}^{*}(\psi _{D_{T}}
(c\theta +D))\frac{\Omega }{\theta
D},\quad \Omega O_{K}). \hskip 2cm  (2.4 )$$
Since $D\equiv 1(mod \  4),\ $ so
$c\theta +D\equiv 1(mod \  \theta )\ $ for any
$c \in \mathcal{C}.\ $  Thus by the definition of
$\psi _{D_{T}}$  and quartic reciprocity we have
$$ \psi _{D_{T}}(c\theta +D)=\overline{\left(\frac{D_{T}}
{c\theta +D} \right) _{4}} (c\theta +D)\\
=\overline{\left(\frac{c\theta +D}{D_{T}}
\right) _{4}}(c\theta +D) \\
=\overline{\left(\frac{c\theta }{D_{T}} \right) _{4}}
(c\theta +D). $$
Then by $(2.4 )$ and the fact
$L_{(\theta D)}(\overline{\psi }_{D_{T}},\quad 1)=
L_{S}(\overline{\psi }_{D_{T}},\quad 1),\ $ we have
$$ \frac{\theta D}{\alpha \omega }L_{S}
(\overline{\psi }_{D_{T}},\quad 1)=
\sum _{c \in \mathcal{C}}E_{1}^{*}
\left( \left(\frac{c\omega }{D}+
\frac{\omega }{\theta } \right)\alpha
\overline{\left(\frac{c\theta }{D_{T}} \right) _{4}},\quad \alpha \omega
O_{K} \right).\hskip 1.6cm (2.5 )$$
Put $\lambda =\alpha \overline{\left(\frac{c\theta }{D_{T}}
 \right) _{4}},\ $ by
$$ E_{1}^{*}(\lambda z,\quad \lambda L)=\lambda ^{-1}E_{1}^{*}(z, \quad L),\quad $$ we
have $$ E_{1} ^{*} ((\frac{c \omega }{D}+\frac{\omega }{\theta }) \alpha
\overline{\left(\frac{c \theta }{D_{T}} \right) _{4}} ,\quad \alpha
\overline{\left(\frac{c \theta }{D_{T}} \right) _{4}} \omega O_{K} )= \frac{1}{\alpha }
\left(\frac{c\theta }{D_{T}} \right) _{4} E_{1} ^{*} (\frac{c\omega }{D}+\frac{\omega
}{\theta },\quad \omega O_{K}) $$ So by $(2.5 )$ we have $$ \frac{\theta D}{\omega }
L_{S}(\overline{\psi } _{D_{T}} ,\quad 1)= \left(\frac{\theta }{D_{T}} \right) _{4} \sum
_{c \in \mathcal{C}} \left(\frac{c}{D_{T}} \right) _{4} E_{1} ^{*} (\frac{c\omega }{D}+
\frac{\omega }{\theta },\quad \omega O_{K}),\qquad (2.6 ) $$ For the period lattice
$L_{\omega }= \omega O_{K}\ $ mentioned above, denote the corresponding Weierstrass $ \wp
-$function by $\wp (z,L_{\omega }),\ $ denote the corresponding Weierstrass $\quad
Zeta-$function by $\zeta (z,L_{\omega }),\ $ then we have $\wp '(z) ^{2}=4\wp (z)
^{3}-4\wp (z) \quad .$ So by results in [Go-Sch]  we have $$ E_{1} ^{*}(\frac{c\omega
}{D}+\frac{\omega }{\theta }, \omega O_{K})\\ = \zeta (\frac{c\omega }{D},L_{\omega
})+\zeta (\frac{\omega }{\theta },L_{\omega })+ \frac{1}{2} \frac{\wp '(\frac{c\omega
}{D})-(2-2I)} {\wp (\frac{c\omega }{D})-I}- \frac{\pi }{\omega }
\overline{(\frac{c}{D}+\frac{1}{\theta })}.\qquad (2.7) $$ The representation system
$\mathcal{C}\ $ of $(O_{K}/(D)) ^{\times }$ may be so chosen that $-c \in \mathcal{C}$
whenever $c \in \mathcal{C}\ $ . Then $\left(\frac{-c}{D_{T}} \right) _{4}=
\left(\frac{c}{D_{T}} \right) _{4}.\ $ Since $\zeta (z,L_{\omega })$  and $\wp
'(z,L_{\omega })$ are odd functions, and $\wp (z,L_{\omega })$
 is even, so by $(2.6 )$ we
have
 $$\frac{D}{\omega }\overline{\left(\frac{\theta }{D_{T}}
 \right) _{4}}
L_{S}(\overline{\psi } _{D_{T}}, 1)
 =\frac{1}{\theta } \left\{\sum _{c \in \mathcal{C}}
\left(\frac{c}{D_{T}} \right) _{4} \zeta (\frac{c\ \omega }{D},
L_{\omega })- \frac{\pi }{\omega }\sum _{c \in \mathcal{C}}
\left(\frac{c}{D_{T}} \right) _{4}
\frac{\overline{c}}{\overline{D}} + \right.$$
$$ \left. \frac{1}{2} \sum _{c \in \mathcal{C}}
\left(\frac{c}{D_{T}} \right) _{4} \frac{\wp '(\frac{c\omega }
{D})} {\wp (\frac{c\ \omega }{D})-I}
- (1-I) \sum _{c \in \mathcal{C}} \left(\frac{c}{D_{T}} \right) _{4}
\frac{1}{\wp (\frac{c\ \omega }{D})-I} \right\}$$
$$ +  \frac{1}{\theta } \sum _{c \in \mathcal{C}}
 \left(\frac{c}{D_{T}} \right) _{4}
\left(\zeta \left(\frac{\omega }{\theta },
L_{\omega }\right)- \frac{\pi }{\omega \overline{\theta
}}\right)$$
$$ =-\frac{1-I}{\theta } \sum _{c \in \mathcal{C}}
 \left(\frac{c}{D_{T}} \right) _{4}
\frac{1}{\wp (\frac{c\ \omega }{D})-I}+
\frac{1}{\theta }\sum _{c \in \mathcal{C}}
\left(\frac{c}{D_{T}} \right) _{4} \left(\zeta
\left(\frac{\omega }{\theta }, L_{\omega }\right)-
\frac{\pi }{\omega \overline{\theta }}\right) $$
That is
$$ \frac{D}{\omega } \overline{\left(\frac{\theta }{D_{T}} \right) _{4}}
L_{S}(\overline{\psi }_{D_{T}},1)=\frac{I}{2} \sum _{c \in \mathcal{C}}
\left(\frac{c}{D_{T}} \right) _{4} \frac{1}{\wp (\frac{c\ \omega }{D})-I}+ $$
$$ \frac{1}{\theta } \sum _{c \in \mathcal{C}} \left(\frac{c}{D_{T}} \right) _{4}
\left(\zeta \left(\frac{\omega }{\theta },L_{\omega }\right)-
\frac{\pi }{\omega \overline{\theta }}\right). \hskip2cm (2.8) $$
By [Zhao]  we know $\zeta (\frac{\omega }{\theta },L_{\omega })-
\frac{\pi }{\omega \overline{\theta }}=\frac{\theta }{4},\ $
so
$$ \frac{D}{\omega } \overline{\left(\frac{\theta }{D_{T}}
\right) _{4}}
L_{S}(\overline{\psi }_{D_{T}},1)=\frac{I}{2}
\sum _{c \in \mathcal{C}}
\left(\frac{c}{D_{T}} \right) _{4}
\frac{1}{\wp (\frac{c\ \omega }{D})-I}+
\frac{1}{4} \sum _{c \in \mathcal{C}}
\left(\frac{c}{D_{T}} \right) _{4}  . $$
This proves Theorem 1.1.
\par \vskip 0.5cm

  {\bf  Lemma  } 2.2 .\quad
$ \sum _{c \in \mathcal{C}} \left(\frac{c}{D_{T}} \right) _{4}
= \left\{ \begin{array}{l} \sharp \mathcal{C}, \qquad  \hbox{if } \quad T=\emptyset  \ ; \\
0, \qquad \hbox{if } \quad T\neq \emptyset  \ .
\end{array}
\right .$
       \par \vskip 0.3cm
{\bf  Proof }$.\quad $ By the definition  of
quartic residue symbol the lemma could be easily verified.
   \par \vskip 0.4cm

{\bf  Lemma   2.3 }.\quad Let
$D=\pi _{1} \cdots \pi _{n}\ $  where $\pi _{k}
\equiv 1(mod \  4)$ are distinct Gaussian prime
$(k=1,\cdots, n) $. Let $c$ be any Gaussian integer
relatively prime to $D$. Then
     \par \vskip 0.2cm
(1) \     $ \sum\limits _{T} \left(\frac{c}{D_{T}} \right) _{4}=
\mu (1+I) ^{t}
\quad \hbox{or } \quad 0  , $
where $\mu \in \{\pm 1,\quad \pm I \}$, $t$ is a integer with
 $n\leq t\leq 2n \ .$  \par \vskip 0.2cm

(2) \     $ \sum\limits _{T} \left(\frac{c}{D_{T}} \right) _{4}
=0 \quad $
 if and only if $\quad \left(\frac{c}{\pi _{k}} \right) _{4}=-1$
 (for some $k \in \{1,\cdots, n \}$).    \par \vskip 0.2cm

(3) \ Suppose that  $\quad \left(\frac{c}{\pi _{k}} \right) _{4}\neq -1$
for arbitrary  $k \in \{1,\cdots, n\}$),then

\hskip 2cm
$ \sum\limits _{T} \left(\frac{c}{D_{T}} \right) _{4}=
\mu (1+I) ^{n+s} , \quad  \quad $
 where  $\mu $ is as in (1) above,

\hskip 2cm  $\quad s=\sharp \{\pi _{k} \ :
\  \ \pi _{k} | D \quad \hbox{and }
\quad \left(\frac{c}{\pi _{k}} \right) _{4}=1,
\quad k=1,\cdots, n \}.$

In particular we know
$$ \sum\limits _{T} \left(\frac{c}{D_{T}} \right) _{4}=2^{n}
\quad \hbox{ if and only if }
\quad \left(\frac{c}{\pi _{1}} \right) _{4} =\cdots =
\left(\frac{c}{\pi _{n}} \right) _{4} =1 \quad ; $$
$$ \sum\limits _{T} \left(\frac{c}{D_{T}} \right) _{4}=
\mu (1+I) ^{n}
\quad \hbox{ if and only if } \quad \left(\frac{c}{\pi _{k}}
\right) _{4} \in
\{I,-I\},\quad k=1,\cdots, n  \ ,   $$
where the sum $\sum\limits _{T}$ is taken over all subsets $T$
of $ \{1,\cdots, n\}$.
         \par \vskip 0.4cm

{\bf  Proof }$.\quad $In fact we have
$\sum\limits _{T} \left(\frac{c}{D_{T}} \right) _{4}=
\left(1+\left(\frac{c}{\pi _{1}} \right) _{4}\right) \cdots
\left(1+\left(\frac{c}{\pi _{n}} \right) _{4}\right),\ $
from which the results could be deduced.
\par \vskip 0.4cm

{\bf Lemma  2.4 }. \quad $v\left( S*(D)\right)\geq
{(n-1)}/{2}$.
\par \vskip 0.4cm

{\bf  Proof }. By results of [Zhao]  or [B-SD],
we know
$$v_2 \left(\wp (\frac{c\omega }{D})-I\right)=\frac{3}{4}$$
(for any Gaussian integer $c  $ relatively prime to
$D$). And by Lemma   2.3   we have
$$ v_2 \left(\sum\limits _{T} \left(\frac{c}{D_{T}}
\right) _{4}\right)=
v_2(\mu (1+I) ^{t})=
\frac{t}{2} \geq \frac{n}{2} \  $$
(Here we regard $v_2(0)$ as $\infty $). Thus by properties of
valuation and our choice of $\mathcal{C}$ with the property
$ c, -c \in \mathcal{C}$,  we have
$$\quad v_2(S^{*}(D))\geq -\frac{3}{4}+\frac{n}{2} . $$
since $\quad \pi _{k}\equiv 1(mod \  4)\quad
(k=1,\cdots, n),\ $ so
$$ N(D_{T})\equiv N(D)\equiv (mod \  8),\quad \left(
\frac{I}{D_{T}} \right) _{4}=I^{(N(D_{T})-1)/4} =\pm 1 \ . $$
Also we have
$$ \sharp (O_{K}/(D)) ^{\times }=\sharp \mathcal{C}=\prod _{k=1} ^{n} (N(\pi
_{k})-1)\equiv 0 (mod \  8)  \ ,  $$
so we could choose $\mathcal{C}\ $ properly such that
$\pm c ,\pm I c \in \mathcal{C}$ (when $c \in \mathcal{C}$).
Put
$$ V=\{c \in \mathcal{C} \ :\  \ c\equiv 1 (mod \ \theta ) \},
\qquad V'=V \cup I V , $$
then $\quad \mathcal{C}=V' \cup (-V')\ .$
Since $\quad I O_{K}=O_{K},\ $ so
$I L_{\omega }=I(\omega O_{K})=\omega O_{K}=L_{\omega }.\ $
 Thus   by the definition of Weierstrass $\wp -$function,
  we could obtain
$$ \wp (Iz,\ I L_{\omega })=\frac{1}{(Iz) ^{2}}+
\sum _{\alpha \in I L_{\omega }}
\left(\frac{1}{(Iz-\alpha ) ^{2}}-
\frac{1}{\alpha ^{2}}\right) $$
$$ = \frac{1}{(Iz) ^{2}}+\sum _{\alpha ' \in L_{\omega }}
\left(\frac{1}{(Iz-I\alpha ')
^{2}}-\frac{1}{(I \alpha ') ^{2}}\right) $$
$$ = \frac{1}{I^{2}} \left(\frac{1}{z^{2}}+
\sum _{\alpha ' \in L_{\omega }}
\left(\frac{1}{(z-\alpha ') ^{2}}-\frac{1}{{\alpha '} ^{2}}
\right) \right)
=-\wp (z,\ L_{\omega }) \ , $$
 that is  $$\quad \wp (Iz,\  L_{\omega })=
 -\wp (z,\ L_{\omega })  .$$
  In particular, assume $\quad z=\frac{c\omega }{D},\ $
 then  we have
$$ \wp \left(\frac{Ic\ \omega }{D},\  L_{\omega }\right)=
-\wp \left(\frac{c\ \omega }{D},
\quad L_{\omega }\right).$$
$$ S^{*}(D)=\frac{I}{2} \sum _{c \in \mathcal{C}} \frac{1}
{\wp (\frac{c\ \omega }{D})-I} \sum _{T} \left(\frac{c}{D_{T}} \right) _{4}=
I \sum _{c \in V'} \frac{1}{\wp (\frac{c\ \omega }{D})-I} \sum _{T}
\left(\frac{c}{D_{T}} \right) _{4} $$
$$ =I \sum _{c \in V} \left [\frac{1}{\wp (\frac{c\ \omega }{D})-I}
\sum _{T} \left(\frac{c}{D_{T}} \right) _{4}+\frac{1}{\wp (\frac{Ic\ \omega }{D})-I}
\sum _{T} \left(\frac{I c}{D_{T}} \right) _{4} \right]  $$
$$ =I \sum _{c \in V} \left [\sum _{T} \frac{1}{\wp (\frac{c\ \omega }{D})-I}
\left(\frac{c}{D_{T}} \right) _{4}+\sum _{T} \frac{1}{-\wp (\frac{c\ \omega }{D})-I}
\left(\frac{I}{D_{T}} \right) _{4} \left(\frac{c}{D_{T}} \right) _{4} \right]  $$
$$ =I\sum _{c \in V} \left [\sum _{T} \left( \frac{1}{\wp (\frac{c\ \omega
}{D})-I}-\left(\frac{I}{D_{T}} \right) _{4} \frac{1}{\wp (\frac{c\ \omega }{D})+I}
\right) \left(\frac{c}{D_{T}} \right) _{4} \right] . $$
Note that $\quad v_2(\wp (({c\ \omega })/{D})-I)={3}/{4},\ $
so we know
$\ v_2(\wp (({c\ \omega })/{D})+I)={3}/{4} \quad .$
Note also
$$ \frac{1}{\wp (\frac{c\ \omega }{D})-I}+\frac{1}{\wp (\frac{c\ \omega
}{D})+I}=\frac{2\wp (\frac{c\ \omega }{D})}{(\wp (\frac{c\ \omega }{D})) ^{2}+
1} \  , $$
$$ \frac{1}{\wp (\frac{c\ \omega }{D})-I}-\frac{1}{\wp (\frac{c\ \omega
}{D})+I}=\frac{2I}{(\wp (\frac{c\ \omega }{D})) ^{2}+1} \ , $$
$$\left(\frac{I}{D_{T}} \right) _{4}=\pm 1 \ , $$
so
$$ S^{*}(D)=I \sum _{c \in V} \frac{2B}{(\wp (\frac{c\ \omega }
{D})) ^{2}+1}
\sum _{T} \left(\frac{c}{D_{T}} \right) _{4} ,\quad B=I,\quad
\hbox{or } \quad
\wp (\frac{c\ \omega }{D}) \ .$$
Since
$$ v_2((\wp (\frac{c\ \omega }{D})) ^{2}+1)=v_2(\wp (\frac{c\ \omega }{D})-I)+
v_2(\wp (\frac{c\ \omega }{D})+I)=\frac{3}{4}+\frac{3}{4}=\frac{3}{2}  \ ,  $$
so
$$ v_2 \left(\frac{2B}{(\wp (\frac{c\ \omega }{D})) ^{2}+1}
\right)=
1-\frac{3}{2}=-\frac{1}{2}  \qquad ( \hbox{and obviously
 we have } v_2(B)=0) \ , $$
 therefore   we have
$$ v_2(S^{*}(D))\geq -\frac{1}{2}+v_2 \left(\sum _{T} \left( \frac{c}{D_{T}} \right) _{4}
\right) \geq -\frac{1}{2}+\frac{n}{2}=\frac{n-1}{2}  . $$
 This proves the lemma.
    \par \vskip 0.5cm

{\bf  Proof of Theorem 1.2 }$.\quad $ First let us prove
$$\quad v_2 \left(L(\overline{\psi }_{D},1)/\omega  \right)
\geq \frac{n-1}{2}  . $$
Take sums for both sides of formula (1.1)
over all subsets $T$ of $\{1,\cdots , n \}$ ,
  we have
$$ \sum _{T} \frac{D}{\omega } \overline{\left(\frac{\theta }{D_{T}} \right) _{4}}
L_{S}(\overline{\psi }_{D_{T}},1)=\frac{I}{2} \sum _{T} \sum _{c \in \mathcal{C}}
\left(\frac{c}{D_{T}} \right) _{4} \frac{1}{\wp (\frac{c\ \omega }{D})-I}+
\frac{1}{4} \sum _{T} \sum _{c \in \mathcal{C}}
\left(\frac{c}{D_{T}} \right) _{4} \quad . $$
 so  By Lemma    2.2   and (1.2) , we obtain
$$ \sum _{T} \frac{D}{\omega } \overline{\left(\frac{\theta }{D_{T}} \right) _{4}}
L_{S} (\overline{\psi } _{D_{T}}, 1)=S^{*}(D)+
\frac{\sharp \mathcal{C}}{4} .\hskip 2.5cm  (2.9)$$
Also we have
$$ v_2(\frac{\sharp \mathcal{C}}{4})=v_2 \left(
\frac{\prod\limits _{k=1}
\limits^{n}(\pi _{k}
\overline{\pi } _{k}-1)}{4} \right) \geq 3n-2\geq n  \ ,  $$
so By Lemma   2.4    we obtain
$$ v_2 \left(\sum _{T} \frac{D}{\omega } \overline{\left(\frac{\theta }
{D_{T}} \right) _{4}} L_{S} (\overline{\psi } _{D_{T}}, 1)\right)
\geq \frac{n-1}{2} \quad. $$
By Lemma   2.1  we know
$L_{S}(\overline{\psi }_{D_{T}},1)=L(\overline{\psi } _{D},1)$
when $T=\{1,\cdots , n\}$; and when
$T=\emptyset $ we have
$$ L_{S}(\overline{\psi } _{D_{T}},1)=L_{S}
(\overline{\psi } _{1},1)=L(\overline{\psi } _{1},1)
\prod\limits _{k=1} \limits^{n}
(1-\frac{\overline{\pi } _{k}}{\pi _{k}\overline{\pi } _{k}})$$
$$ =L(\overline{\psi } _{1},1) \prod \limits _{k=1}^{n}
(1-\frac{1}{\pi _{k}}) .$$
By [B-SD] or [Zhao]  we know
 $\quad L(\overline{\psi }_{1},1)={ \omega }/{4},\ $
 so  $$\quad L_{S}(\overline{\psi } _{1},1)=
 \frac{\omega }{4} \prod \limits
_{k=1} ^{n}(1-\frac{1}{\pi _{k}}),\ $$
$$ v_2 \left(L_{S}(\overline{\psi } _{1},1)/\omega \right)=
v_2 \left(\frac{1}{4} \prod _{k=1} ^{n} (1-\frac{1}{\pi _{k}})
 \right)\geq 2n-2 \qquad ( \hbox{since }
\quad v_2(\pi _{k}-1)\geq 2 ). $$
Now we use induction on $n$ to prove our assertion
$ v_2 \left(L(\overline{\psi }_{D},1)/\omega  \right)\geq
\frac{n-1}{2}$.
 If $n=1$, then $D=\pi _{1},\ $
 $L_{S}(\overline{\psi } _{1},1)=({\omega }/{4})\cdot
({\pi _{1}-1})/{\pi _{1}}.\ $
Since $\pi _{1}\equiv 1 (mod \  4),\ $ so
$v_2 \left(L_{S}(\overline{\psi } _{1},1)/\omega \right)
\geq 0 \quad .$ By the above analysis we have
$$ v_2 \left(\frac{\pi _{1}}{\omega } \left(\frac{\theta }{1} \right) _{4}L_{S}
(\overline{\psi } _{1},1)+\frac{\pi _{1}}{\omega }\overline{\left(\frac{\theta }
{\pi _{1}} \right) _{4}} L_{S}(\overline{\psi }_{\pi _{1}},1) \right)\geq
\frac{1-1}{2}=0  \ ,  $$
 therefore   we have
$$ v_2 \left(L(\overline{\psi } _{\pi _{1}},1)/\omega \right)=
v_2 \left(L_{S}(\overline{\psi } _{\pi _{1}},1)/
\omega \right)\geq 0 \quad . $$
Now assume our assertion
 is true for cases
$1,2,\cdots , n-1$, and consider the case $n$,
 $\quad D=\pi _{1} \cdots \pi _{n}.\ $
 For any subset $T$ of $\{1,\cdots , n\}$,
denote $t=t(T)=\sharp T$, by Lemma  2.1  we know
$$\frac{D}{\omega } \overline{\left(\frac{\theta }{D_{T}} \right) _{4}}
L_{S} (\overline{\psi } _{D_{T}}, 1)=\frac{D}{\omega }
\overline{\left(\frac{\theta }{D_{T}} \right) _{4}}
L (\overline{\psi } _{D_{T}}, 1) \prod _{\pi _{k}|\widehat{D}_{T}}
\left(1-\left(\frac{D_{T}}{\pi _{k}} \right) _{4}
\frac{1}{\pi _{k}} \right)  $$
Since $\left({D_{T}}/{\pi _{k}} \right) _{4}=
\pm 1 ,\pm I,\ $ so
$$ 1-\left(\frac{D_{T}}{\pi _{k}} \right) _{4} \frac{1}{\pi _{k}}=
\frac{\pi _{k}-\mu }{\pi _{k}}  \qquad  \mu \in \{\pm 1 ,\pm I \}. $$
Note that  $\pi _{k}\equiv 1 (mod \ 4),\ $ so  we know
 $\quad v_2(\pi _{k}-\mu )\geq
{1}/{2};\ $ moreover the equality holds
 if and only if $\left({D_{T}}/{\pi _{k}} \right) _{4} ^{2}
=-1 \quad .$ Thus by our inductive assumption, we have
$$ v_2 \left(\frac{D}{\omega } \overline{\left(
\frac{\theta }{D_{T}} \right) _{4}}
L_{S} (\overline{\psi } _{D_{T}}, 1) \right)=
v_2 \left(L(\overline{\psi }
_{D_{T}},1)/\omega  \right)+
\sum _{\pi _{k}|\widehat{D}_{T}}
v_2 \left(1-\left(\frac{D_{T}}{\pi _{k}} \right) _{4}
\frac{1}{\pi _{k}}  \right)$$
$$ \geq \frac{t-1}{2}+\frac{1}{2}\cdot
\sharp \{ \pi _{k} \ :\  \ \pi _{k}|\widehat{D}_{T} \}=
\frac{t-1}{2}+\frac{n-t}{2}=\frac{n-1}{2} . $$
Also  when $T=\emptyset $ we have
$$ L_{S}(\overline{\psi } _{D_{T}},1)=
L_{S}(\overline{\psi } _{1},1)=
L(\overline{\psi } _{1},1) \prod _{k=1} ^{n}
(1-\frac{1}{\pi _{k}}) =
\frac{\omega }{4} \prod _{k=1} ^{n} (1-\frac{1}{\pi _{k}}), $$
therefore
$$ v_2 \left( L_{S} (\overline{\psi } _{1},1)/\omega \right)\geq 2n-2\geq
\frac{n-1}{2}  \ ,  $$
\begin{align*}
&v_2 \left(L(\overline{\psi }_{D},1)/\omega  \right)=
v_2 \left(\frac{D}{\omega } \overline{\left(
\frac{\theta }{D} \right) _{4}}
L (\overline{\psi } _{D}, 1) \right) \\
& \qquad =v_2 \left(\sum _{T} \frac{D}{\omega }
\overline{\left(\frac{\theta }{D_{T}} \right) _{4}}
L_{S} (\overline{\psi } _{D_{T}}, 1)-\sum _{\emptyset
\neq T \subsetneqq \{1,\cdots , n\}} \frac{D}{\omega }
\overline{\left(\frac{\theta }{D _{T}} \right) _{4}}
L_{S} (\overline{\psi } _{D_{T}}, 1) -\frac{D}{\omega }
L_{S}(\overline{\psi }
_{1},1) \right)\\
& \qquad \geq \frac{n-1}{2}  .
\end{align*}
 Therefore  by induction we have proved our
 assertion  that
$ v_2 \left(L(\overline{\psi }_{D},1)/\omega  \right)\geq
\frac{n-1}{2}$ holds for any positive integer $n$.
   \par \vskip 0.3cm

Now we consider the condition for the equality holds,
using induction method on $n$ too.
 If $n=1$, then $D=\pi _{1},\ $ by (2.9) we obtain
$$ \frac{\pi _{1}}{\omega } \overline{\left(\frac{\theta }{1} \right) _{4}}
L_{\pi _{1}}(\overline{\psi } _{1},1)+\frac{\pi _{1}}{\omega }
\overline{\left(\frac{\theta }{\pi _{1}} \right) _{4}}
L_{\pi _{1}}(\overline{\psi } _{\pi _{1}},1)=S^{*}(\pi _{1})+\frac{\pi _{1}
\overline{\pi _{1}}-1}{4} \ , $$
 that is
$$ \qquad \frac{1}{4} (\pi _{1}-1)+\frac{\pi _{1}}{\omega }
\left(\frac{\theta }{\pi _{1}} \right) _{4} L(\overline{\psi } _{\pi _{1}},1)
=S^{*}(\pi _{1})+\frac{\pi _{1} \overline{\pi _{1}}-1}{4} \ .$$
Since
$$ v_2 \left(\frac{\pi _{1} \overline{\pi _{1}}-1}{4} \right)=
v_2(\pi _{1} \overline{\pi _{1}}-1)-2\geq 1 ,\quad v_2(S^{*}(\pi _{1}))\geq
\frac{1-1}{2}=0 \quad (\hbox{Lemma  } 2.4 ) , $$
so the equality
\begin{align*}
v_2 \left(L(\overline{\psi }_{\pi _{1}},1)/\omega  \right)&=
v_2 \left(\frac{\pi _{1}}{\omega } \overline{\left(\frac{\theta }{\pi _{1}}
\right) _{4}} L _{\pi _{1}} (\overline{\psi } _{\pi _{1}}, 1) \right)\\
&=
v_2 \left(S^{*}(\pi _{1})+\frac{\pi _{1} \overline{\pi _{1}}-1}{4} -
\frac{1}{4} (\pi _{1}-1) \right)\\
&=0
\end{align*}
holds if and only if one of the following conditions is true:

\quad (1) $\ v_2(\pi _{1}-1)=2 \quad $
when $\ v_2(S^{*}(\pi _{1}))>0 $;

\quad (2) $\ v_2(\pi _{1}-1)>2 \quad $
when $\ v_2(S^{*}(\pi _{1}))=0 $.

Thus we know
 $$ v_2 \left(L(\overline{\psi }_{\pi _{1}},1)/
 \omega  \right)=0
 \hbox{ holds if and only if }
 \delta _{1}(\pi _{1})=s_{1}(\pi _{1})+
\varepsilon _{1}(\pi _{1})=1.$$

Assume our result is true in the cases $1,\cdots , n-1$,
consider the case $n,\ $ i.e.
$\ D=\pi _{1} \cdots \pi _{n} \ .$
When $T=\emptyset ,\ $ we have
$$ \frac{D}{\omega } \overline{\left(\frac{\theta }{D_{T}} \right) _{4}}
L_{S} (\overline{\psi } _{D_{T}}, 1)=\frac{D}{\omega }
\overline{\left(\frac{\theta }{D_{T}} \right) _{4}}
L_{D} (\overline{\psi } _{1}, 1)= \frac{D}{\omega } L (\overline{\psi } _{1},1)
\prod _{k=1} ^{n}(1-\frac{1}{\pi _{k}}) , $$
\begin{align*} v_2 \left(\frac{D}{\omega } \overline{\left(\frac{\theta }
{D_{T}} \right) _{4}}
L_{S} (\overline{\psi } _{D_{T}}, 1)\right)&= v_2 \left(L(\overline{\psi }
_{1},1)/\omega  \right)+\sum _{k=1} ^{n} v_2(\pi _{k}-1)\\
&=v_2(\frac{1}{4})+
\sum _{k=1} ^{n} v_2(\pi _{k}-1)\\
&\geq 2n-2 \geq n \geq \frac{n-1}{2} \ .
\end{align*}
When $\quad \emptyset \neq T \subsetneqq \{1,\cdots , n\}$
we have ,
\begin{align*} \hskip 3cm v_2 \left(\frac{D}{\omega }
 \overline{\left(\frac{\theta }{D_{T}} \right) _{4}}
L_{S} (\overline{\psi } _{D_{T}}, 1)\right)&
 = v_2 \left(L_{S}(\overline{\psi }
_{D_{T}},1)/\omega  \right)\\
&\hskip -3cm =v_2 \left(\frac{L(\overline{\psi } _{D_{T}},1)}
{\omega } \cdot \prod \limits _{\pi _{k}|\widehat{D}_{T}}
\left(1-\left(\frac{D_{T}}{\pi _{k}} \right) _{4} \frac{1}
{\pi _{k}} \right) \right)\\
&\hskip -3cm= v_2 \left(L(\overline{\psi } _{D_{T}},1)
/\omega  \right)+\sum \limits
_{\pi _{k}|\widehat{D} _{T}} v_2 \left(1-\left(\frac{D_{T}}
{\pi _{k}} \right) _{4}
\frac{1}{\pi _{k}} \right) .
\end{align*}
Since $\ \left({D_{T}}/{\pi _{k}} \right) _{4}=\pm 1,\ $
 $\pm I,\ $ so
$$ 1- \left(\frac{D_{T}}{\pi _{k}} \right) _{4} \frac{1}{\pi _{k}}=
\frac{\pi _{k}-\mu }{\pi _{k}}=\frac{\pi _{k}-1+(1-\mu )}{\pi _{k}} ,\quad
\mu \in \{\pm 1 , \pm I \}.$$
 Therefore  $\quad v_2 \left(1- \left(\frac{D_{T}}{\pi _{k}} \right) _{4} \frac{1}
{\pi _{k}} \right)\geq \frac{1}{2},\ $
 and the equality holds if and only if
$\quad \left({D_{T}}/{\pi _{k}} \right) _{4} =\pm I ,\ $ \\
 i.e. $\quad \left({D_{T}}/{\pi _{k}} \right) _{4} ^{2} =-1,\ $
 that is
$\quad \left[\frac{D_{T}}{\pi _{k}} \right]  _{2}=1  .$
Thus
$$\quad v_2 \left(1- \left(\frac{D_{T}}{\pi _{k}} \right) _{4} \frac{1}
{\pi _{k}} \right)=\frac{1}{2} \quad \hbox{ if and only if }
\quad  \left[\frac{D_{T}}{\pi _{k}} \right]  _{2}=1  . $$
By the proof of the first part of the Theorem we know
$$ v_2 \left(L(\overline{\psi } _{D_{T}},1)/\omega  \right)\geq
\frac{t(T)-1}{2},\quad t(T)=\sharp T  \ ,  $$
and by our inductive assumption we know the equality holds
 if and only if $\quad \delta _{t} (D_{T})=1 ,
 \quad t=t(T) .$
Thus
$$v_2 \left(\frac{D}{\omega } \overline{\left(\frac{\theta }
{D_{T}} \right) _{4}} L_{S} (\overline{\psi } _{D_{T}}, 1)\right)\geq
\frac{t(T)-1}{2}+\frac{n-t(T)}{2}=\frac{n-1}{2} , $$
and the equality holds if and only if
$\quad  \left[\frac{D_{T}}{\pi _{k}} \right]  _{2}=1 \quad $
( for any $\ \pi _{k}|\widehat{D}_{T})\ $ and
$\ \delta _{t}(D_{T})=1 \ .$
That is to say the equality
$$ v_2 \left(\frac{D}{\omega } \overline{\left(\frac{\theta }
{D_{T}} \right) _{4}} L_{S} (\overline{\psi } _{D_{T}}, 1)
\right)=\frac{n-1}{2} $$ holds  if and only if
$$  \left(\prod _{\pi _{k}|\widehat{D} _{T}}
\left[\frac{D_{T}}{\pi _{k}} \right]  _{2} \right)
\delta _{t}(D_{T})=1 .$$
For the  elliptic curve $E_{D_{T}}:\ y^{2}=x^{3}-D_{T}x \quad $
  and Hecke characters $\psi _{D_{T}},\ $
   by  [Ru\ 1-2]  we know $L(\overline{\psi }_{D_{T}},1)/\Omega
\in K=\Q(I),\ $ also we have $ \quad \Omega =
\frac{\omega }{\sqrt[4] {D_{T}}}  \ ,  $ so
\begin{align*} L(\overline{\psi }_{D_{T}},1)/ \omega
&=(\sqrt[4] {D_{T}}) ^{-1} \cdot
L(\overline{\psi }_{D_{T}},1)/ \frac{\omega }{\sqrt[4] {D_{T}}}\\
&=(\sqrt[4] {D_{T}}) ^{-1} \cdot L(\overline{\psi }_{D_{T}},1)/\Omega \in
K(\sqrt[4] {D_{T}}),
\end{align*}
 i.e. $\quad  L(\overline{\psi }_{D_{T}},1)/ \omega \in
K(\sqrt[4] {D_{T}})  . \ $
Thus by Lemma   2.1    we know
$$ \frac{D}{\omega } \overline{\left(\frac{\theta }
{D_{T}} \right) _{4}} L_{S} (\overline{\psi } _{D_{T}}, 1)=
D \overline{\left(\frac{\theta }{D_{T}} \right) _{4}}
\prod _{\pi _{k}|\widehat{D} _{T}} \left(1-\left(\frac{D_{T}}{\pi _{K}}
\right) _{4} \frac{1}{\pi _{k}} \right) \cdot L(\overline{\psi }_{D_{T}},1)/
\omega \in K(\sqrt[4] {D_{T}}) , $$
and if  $\quad  v_2 \left(\frac{D}{\omega } \overline{\left(\frac{\theta }
{D_{T}} \right) _{4}} L_{S} (\overline{\psi } _{D_{T}}, 1)\right)=
\frac{n-1}{2},\ $ then we have
$$ \frac{D}{\omega } \overline{\left(\frac{\theta }
{D_{T}} \right) _{4}} L_{S} (\overline{\psi } _{D_{T}}, 1)=
(1+I) ^{n-1} \alpha _{_{T}} \sqrt[4] {D_{T}^{3}}  ,  $$
where $\ \alpha _{_{T}} \in K,\ $ and
$\ v_2(\alpha _{_{T}})=0 .$
(since $\quad v_2(\sqrt[4] {D_{T}} ^{3})=
\frac{3}{4} v_2(D_{T})=0$). For any subsets $T$ and $T'\ $ of
$\{1,\cdots , n \}$,  if $\ v_2(\alpha _{_{T}})=
\ v_2(\alpha _{_{T'}})=0 ,\ $ then it could be easily
 verified  that
$$\quad v_2 \left(\alpha _{_{T}} \sqrt[4] {D_{T} ^{3}}+
\alpha _{_{T'}} \sqrt[4] {D_{T'} ^{3}} \right) > 0 \quad . $$
 Thus, consider the terms in the sum
 $$ \sum _{\emptyset \neq T \subsetneqq \{1,\cdots , n\}}
\frac{D}{\omega } \overline{\left(\frac{\theta }
{D_{T}} \right) _{4}} L_{S} (\overline{\psi } _{D_{T}}, 1), $$
for any two terms  with 2-adic valuations  equal to
$(n-1)/2$,
the 2-adic valuation of their sum would be bigger than
  $(n-1)/2$.  So when $n>1$ we have
$$ v_2 \left(\frac{D}{\omega } \overline{\left(\frac{\theta }
{D_{\emptyset }} \right) _{4}} L_{S}
(\overline{\psi } _{D_{\emptyset }}, 1)\right)\geq 2n-2
\geq n>\frac{n-1}{2}  .$$
Hence we know
$v_2 \left(L(\overline{\psi } _{D},1)/\omega \right)=
\frac{n-1}{2} \ $ holds  if and only if
one of the following statements is true :

$(1)\ v_2(S^{*}(D))> ({n-1})/{2}  ;\ $

$(2)\ v_2(S^{*}(D))= \frac{n-1}{2} .\ $

 Statement (1) means that, when
   $\quad \varepsilon _{n}(D)=0$,
in the sum
$$\quad \sum _{\emptyset \neq T \subsetneqq \{1,\cdots , n\}}
\frac{D}{\omega } \overline{\left(\frac{\theta }
{D_{T}} \right) _{4}} L_{S} (\overline{\psi } _{D_{T}}, 1), $$
the number of terms with 2-adic valuation $(n-1)/2$ is odd.
The number of such terms turns to be
\begin{align*}  &\sharp \left\{\emptyset \neq T \subsetneqq \{1,\cdots , n \} :\
v_2 \left(\frac{D}{\omega } \overline{\left(\frac{\theta }
{D_{T}} \right) _{4}} L_{S} (\overline{\psi } _{D_{T}},
1)\right)=\frac{n-1}{2} \right\}\\
&=\sharp \left\{\emptyset \neq
T \subsetneqq \{1,\cdots , n \} :\ \left(\prod _{\pi _{k}|\widehat{D} _{T}}
\left[\frac{D_{T}}{\pi _{k}} \right]  _{2} \right) \delta _{t}(D_{T})=1
\right\}\\
&\equiv \sum _{\emptyset \neq T \subsetneqq \{1,\cdots , n \} }
\left(\prod _{\pi _{k}|\widehat{D} _{T}}
\left[\frac{D_{T}}{\pi _{k}} \right]  _{2} \right) \delta _{t}(D_{T})\\
&\equiv 1(mod\ 2),
\end{align*}
we have
$$ \delta _{n}(D)=\varepsilon _{n}(D)+\sum _{\emptyset
\neq T \subsetneqq \{1,\cdots , n \} }
\left(\prod _{\pi _{k}|\widehat{D} _{T}}
\left[\frac{D_{T}}{\pi _{k}} \right]  _{2} \right) \delta _{t}(D_{T})
\equiv 1(mod\ 2)  . $$

On the other hand, the statement (2) above means that,
  when  $\varepsilon _{n}(D)=1$, in the sum
$$\quad \sum _{\emptyset \neq T \subsetneqq \{1,\cdots , n\}}
\frac{D}{\omega } \overline{\left(\frac{\theta }
{D_{T}} \right) _{4}} L_{S} (\overline{\psi } _{D_{T}}, 1), $$
the number of terms with 2-adic valuation $(n-1)/2$ is even,
that is
$$ \sum _{\emptyset \neq T \subsetneqq \{1,\cdots , n \} }
\left(\prod _{\pi _{k}|\widehat{D} _{T}}
\left[\frac{D_{T}}{\pi _{k}} \right]  _{2} \right)
\delta _{t}(D_{T}) \equiv 0(mod\ 2) .$$
We have
$$ \delta _{n}(D)=\varepsilon _{n}(D)+\sum _{\emptyset
\neq T \subsetneqq \{1,\cdots , n\} }
\left(\prod _{\pi _{k}|\widehat{D} _{T}}
\left[\frac{D_{T}}{\pi _{k}} \right]  _{2} \right) \delta _{t}(D_{T})
\equiv 1(mod\ 2)  . $$
By the discussion above we know that
$\quad v_2 \left(L(\overline{\psi } _{D},1)/\omega  \right)=
({n-1})/{2}
\ $ if and only if $\ \delta _{n}(D)=1 \ .$
 This proves the theorem.

 \par \vskip 0.5cm

{\bf  Proof of Theorem 1.3}. \quad
This theorem follows from  above Theorem 1.2
 and the main result of Coates-Wiles
in [Co-Wi] .
 \par \vskip 0.5cm

 For elliptic curve $\quad E_{D}:\ y^{2}=x^{3}-D x ,\ $
with $\quad D=\pi _{1} ^{2} \cdots \pi _{r} ^{2} \pi _{r+1}
 \cdots \pi _{n},\ $
where
$\ \pi _{k}\equiv 1 (mod\ 4) \ $
are distinct Gaussian prime integers ( $k=1,\ \cdots ,\ n),\ $
we could prove Theorem 1.4 and Theorem 1.5 similarly as
proving Theorem 1.1 and 1.2.
\par \vskip 1cm

\begin{center}
{\Large {\bf   III. 3-Valuations of $L-$series of Elliptic
Curves with CM by $\sqrt {-3}$ }}\\
\end{center}

Now assume the number field  $  K=\Q(\sqrt{-3}),\ $
$\tau =({-1+\sqrt{-3}})/{2}=exp({2\pi I/3})$ is a
 primitive  cubic root of unity, $O_{K}=\Z[\tau ] $ is
the ring of integers of $K$, $I=\sqrt{-1} \quad .$
We now study elliptic curves with complex multiplication
by $ \sqrt{-3}\ $and prove Theorem 1.6 and 1.7.
    \par \vskip 0.3cm

Consider elliptic curve
 $\quad E=E_{D^2}:\ y^{2}=x^{3}+D ,\ $ $D \in O_{K} . $
Let $\psi _{E/K}$ be the Hecke character of $K=\Q(\sqrt{-3}) $
corresponding to the elliptic curve $E/K$
(The discriminant of $E$ is $\Delta  (E)=-2^{4}3^{3}D^{2}$).
For any prime ideal $\wp \dag 6D\ $ of $K$, put
$${ \cal A} _{\wp }=N _{K/\Q} (\wp )+1-\sharp
\widetilde{E}(\F_{\wp }),\ $$
 where $\F_{\wp }=O_{K}/\wp $ is the residue field of
  $K$ modulo $\wp $,  the symbol `` $N _{K/\Q}(\cdot )$''
( or `` $N(\cdot )$'' ) denote the norm map of ideals from $K$
to $\Q$, $\ \widetilde{E}$ is the reduction curve of the
elliptic curve $E$ modulo $\wp $. Then by [Sil\ 2] and
 [Sil\ 1] , we have
$$\quad N _{K/\Q} (\wp )=N _{K/\Q}(\psi _{E/K}(\wp )), $$
$$\sharp \widetilde{E}(\F_{\wp })= N _{K/\Q} (\wp )+
1-\psi _{E/K}(\wp )-\overline{\psi _{E/K}(\wp )} , $$
$$ \quad {\cal A} _{\wp }=\psi _{E/K}(\wp )+
\overline{\psi _{E/K}(\wp )} .$$
By definition, we know that the $L-$series $E_{D}$
(omitting Euler factors corresponding to bad reductions)is
\begin{align*}  L_{D}(s)&=\prod _{\wp \dag 6D} (1-{\cal A} _{\wp } N (\wp )^{-s}+
N(\wp )^{1-2s}) ^{-1}\\
&=\prod _{\wp \dag 6D} (1+(\sharp \widetilde{E}(\F_{\wp })-N(\wp )-1)
N(\wp )^{-s}+N(\wp )^{1-2s}) ^{-1} \\
&=\prod _{\wp \dag 6D}(1-(\psi _{E/K}(\wp )+\overline{\psi _{E/K}(\wp )})
N(\wp )^{-s}+N(\wp )^{1-2s}) ^{-1} \\
&=\prod _{\wp \dag 6D}(1-(\psi _{E/K}(\wp )N(\wp )^{-s}) ^{-1} \cdot
\prod _{\wp \dag 6D}(1-\overline{\psi _{E/K}(\wp )} N(\wp )^{-s}) ^{-1}
\end{align*}
Also by   [Sil\ 2]  we know , $$\psi _{E/K}(\wp )=\overline{\left(\frac{4D}{\pi }
\right) _{6}} \pi ,\ $$
 where $\wp =(\pi )$and $\pi \equiv 1 (mod\ 3)$
 is a primitive prime element . So
\begin{align*}  L_{D}(s)&=\prod _{\hbox{prime }\ \pi \equiv 1(mod\ 3)}
\left(1-\left(\overline{\left(\frac{4D}{\pi }\right) _{6}} \pi +
\left(\frac{4D}{\pi }\right) _{6} \overline{\pi }\right) N(\pi )^{-s}+
N(\pi )^{1-2s} \right) ^{-1} \\
&=\prod _{\pi \equiv 1(mod\ 3)} \left[\left(1-\overline
{\left(\frac{4D}{\pi }\right) _{6}} \cdot \frac{\pi }{(N(\pi )) ^{s}}
\right) ^{-1} \cdot \left(1-\left(\frac{4D}{\pi }\right) _{6}
\cdot \frac{\overline{\pi }}{(N(\pi )) ^{s}} \right) ^{-1} \right]  \\
&=L(s,\psi _{E/K}) \cdot L(s,\overline{\psi }_{E/K}) \quad .
\end{align*}
 Thus, ignoring Euler factors corresponding to primes
 $\wp |6D$,  we have
 $$ \qquad L_{D}(s)=L(s,\psi _{E/K}) \cdot
 L(s,\overline{\psi }_{E/K})
. \hskip 4cm (3.1)$$

Now assume  $D$ is  a  rational  integer ,
Consider the elliptic curve
$\quad E_{D}:\ y^{2}=x^{3}+D,\ $ $D \in \Z . $
For any    rational prime $\quad p\dag 6D ,\ $
$E_{D}$ has good reduction at $p$ ( Since the
discriminant of $E_{D}$ is $-2^{4}3^{3}D^{2}$).
Let $\quad N _{p}=\sharp \widetilde{E}_{D}(\F_{p}),\ $
 i.e. the number of $\F_{p}-$rational points
of the reduced curve $\widetilde{E}_{D}$  . We know

 (1)\ when $\ p\equiv 2(mod\ 3)$, then $N_{p}=p+1 \quad ; $

(2)\  when $\quad p\equiv 1(mod\ 3)$, then
$N_{p}=p+1-\left(\frac{4D} {\overline{\pi }} \right) _{6}
\pi -\left(\frac{4D} {\pi } \right) _{6} \overline{\pi },\ $

where $p=\pi \overline{\pi }=N_{K/\Q}(\pi ),\ $
$\pi \in \Z[\tau ]  \ $ and $\ \pi \equiv 1(mod\ 3)$
is a prime element
( see [Ire-Ro] ). Thus
\begin{align*}  L_{D}(s)
&=\prod _{p\equiv 2(mod\ 3)} (1+p^{1-2s}) ^{-1}\cdot
\prod _{p\equiv 1 (mod\ 3)} \left(1-\pi \left(\frac{4D}
{\overline{\pi }} \right) _{6}p^{-s}-\overline{\pi }\left(
\frac{4D}{\pi } \right) _{6}p^{-s}+p^{1-2s} \right) ^{-1} \\
&=\prod _{\hbox{prime }\ \pi \equiv 1 (mod\ 3)}
\left(1-\left(\frac{4D}{\pi }
\right) _{6}\cdot \frac{\overline{\pi }}{(N(\pi )) ^{s}}
\right) ^{-1}
\quad .
\end{align*}
where  the product  $\prod $ is taken over primitive
prime elements prime $\pi (\equiv 1 (mod\ 3))$ of
 $K=\Q(\sqrt{-3})$.  So
$$ L_{D}(s)=\sum _{\sigma \equiv 1(mod\ 3)} \left(\frac{4D}
{\sigma }\right) _{6}\cdot \frac{\overline{\sigma }}
{(N(\sigma )) ^{s}}
\quad . \hskip 3cm  (3.2) $$
where the sum  $\Sigma $ is taken over primitive integers
 $\sigma (\equiv 1(mod\ 3))$ in $ K$.

{\bf   Remark }  3.1 .
The $L_{D}(s)$ in  $(3.1)$ denotes the
$L-$series  of $E/K$ ,  i.e. $\ L_{D}(s)=L(E/K, s).\ $
While the $L_{D}(s)$ in  $(3.2)$ denotes the
$L-$series of $E/\Q$ ,  i.e. $\  L_{D}(s)=L(E/\Q, s) \ .$
So if $D$ is a  rational integer, by these two  formulae
we obtain $\quad L(E/K, s)=(L(E/\Q, s))^{2}$ (a square), and
$$L(E/\Q, s)=L(s,\psi _{E/K})=L(s,\overline{\psi }_{E/K}) $$
(ignoring Euler factors corresponding to prime$\wp |6D$).
  \par \vskip 0.4cm

Now assume $D \in O_{K}=\Z[\tau ] ,\ $  and $D\equiv 1(mod\ 3)$
 is a primitive integer , $\Delta  $ is the square-free part
 of $D$, i.e. the product of all distinct prime divisors of
$D$, and $\ \Delta  \equiv 1(mod\ 3).\ $
Assume $L_{D}(s)$ as in $(3.2)$. Then For elliptic curve
$\quad E_{D}:\ y^{2}=x^{3}+ 2^{4}D,\ $  we have
$$ L_{D}(s)=\sum _{\sigma \equiv 1(mod\ 3)} \left(\frac{D}
{\sigma } \right) _{6}\cdot \frac{\overline{\sigma }}
{(N(\sigma )) ^{s}} . $$
Since $\quad (3,\Delta  )=1,\ $ so
$\sigma =v\cdot 3\Delta  +3\beta +\Delta ,\ $
where $v$ is an algebraic  integer in $K=\Q(\sqrt{-3})$,
$\beta $ runs over a residue system of $O_K$ modulo $\Delta $.
By the cubic reciprocity  we obtain
$$ \left(\frac{D}{\sigma } \right) _{6} ^{2}=\left(\frac{D}{\sigma } \right) _{3}
= \left(\frac{\sigma }{D } \right) _{3}=
\left(\frac{v\cdot 3\Delta  +3\beta +\Delta  }{D } \right) _{3}=
\left(\frac{3\beta }{D } \right) _{3}=
\left(\frac{3\beta }{D } \right) _{6} ^{2}  \ ,  $$
so $\quad \left(\frac{D}{\sigma } \right) _{6}=
\pm \left(\frac{3\beta }{D } \right) _{6},\ $  that means
the value of $\ \left(\frac{D}{\sigma } \right) _{6}$
depends only on $\beta, \ $ not on $v$.  Thus
\begin{align*}  L_{D}(s)&= \sum _{\sigma \equiv 1(mod\ 3)} \left(\frac{D}{\sigma }
\right) _{6}\cdot \frac{\overline{\sigma }}{(N(\sigma )) ^{s}}\\
&=\sum _{\sigma \equiv 1(mod\ 3)} \left(\frac{D}{\sigma }
\right) _{6}\cdot \frac{3\overline{v}\overline{\Delta  }+3\overline{\beta }+
\overline{\Delta  }}{N(3v\Delta  +3\beta +\Delta  )^{s}}\\
&=\sum _{\beta } \left(\frac{D}{\sigma }
\right) _{6}\sum _{v} \frac{3\overline{v}\overline{\Delta  }+3\overline{\beta }+
\overline{\Delta  }}{N(3v\Delta  +3\beta +\Delta  )^{s}}\\
&=\sum _{\beta } \left(\pm \left(\frac{3\beta }{D} \right) _{6} \right)
\sum _{v} \frac{3\overline{v}\overline{\Delta  }+3\overline{\beta }+
\overline{\Delta  }}{N(3v\Delta  +3\beta +\Delta  )^{s}} .
\end{align*}
 So we have
$$ L_{D}(s)=\sum _{\beta } \left(\pm \left(\frac{3\beta }{D} \right) _{6} \right)
\sum _{v} \frac{3\overline{v}\overline{\Delta  }+3\overline{\beta }+
\overline{\Delta  }}{N(3v\Delta  +3\beta +\Delta  )^{s}} \hskip 2cm  (3.3 ) $$
The function $L_{D}(s)$ defined by these formula is
convergent  for $Re(s)>\frac{3}{2}$.
\par \vskip 0.4cm

Now consider the analytic extension of the above
$L-$series $L_{D}(s)$. we could analytically extend it near
to $s=1$, and express it as a finite sum.
First, for the period lattice $L=O_{K}=\Z+\Z \tau $
 and $z \in \C -L,\ $ Define the following series
 as in [Go-Sch] :
$$ \psi (z, s, L)=\frac{\overline{z}}{|z| ^{2s}}+\sum _{\alpha  \in L - \{0\}}
\left\{ \frac{\overline{z}+\overline{\alpha }}{|z+\alpha | ^{2s}}-
\frac{\overline{\alpha }}{|\alpha |^{2s}} \left(1-\frac{sz}{\alpha }-\frac{(s-1)\overline{z}}
{\overline{\alpha }} \right) \right\}  \qquad (3.4 )$$
 This series is convergent and defines an analytic
 function for $\quad Re(s)>\frac{1}{2}\quad $,
  , and  $\psi (z, s, L)$ is  uniformly convergent to
 $\zeta (z, L)\ $ (the Weierstrass $  \quad Zeta$function
 over the period lattice $L$) when $s\rightarrow 1$ , that is
\begin{align*}  \psi (z,1,L)&=\frac{\overline{z}}{|z| ^{2}}+\sum _
{\alpha \in L-\{0\}} \left\{\frac{\overline{z}+\overline{\alpha }}{|z+\alpha |
^{2}}-\frac{\overline{\alpha }}{|\alpha| ^{2}} (1-\frac{z}{\alpha})
\right\} \\
&=\frac{1}{z}+\sum _{\alpha \in L-\{0\}} \left\{ \frac{1}{z+\alpha }-\frac{1}
{\alpha }+\frac{z}{\alpha ^{2}} \right\} \\
&=\zeta (z,L).
\end{align*}

Also when $Re(s)$ is sufficiently large, we know that
$$ \sum \limits _{\alpha \in L-\{0\}} \frac{\overline{\alpha }}
{|\alpha | ^{2s}}\hskip 1.5cm
 \hbox{and} \hskip 1.5cm
 \sum \limits _{\alpha \in L-\{0\}} \frac{1}{|\alpha | ^{2s}}
  \cdot  \frac{\overline{\alpha }}{\alpha}$$
are absolutely convergent,  so their terms may be re-arranged.
Now we show that when $Re(s)$ is sufficiently large we have
$$ \sum _{\alpha \in L-\{0\}} \frac{\overline{\alpha }}{|\alpha | ^{2s}}
= \sum _{\alpha \in L-\{0\}} \frac{1}{|\alpha | ^{2s}} \cdot \frac{\overline
{\alpha }}{\alpha}=0 \qquad (3.5)$$
In fact, we know that the unit group of $O_{K}$ is
 $\{\pm 1,\pm \tau ,\pm \tau ^{2}\},\ $ and the series is
 absolutely convergent, so in the summation over
the integers $\alpha$,
we could first add all the terms corresponding to
the associates
$\pm \alpha ,\pm\tau \alpha ,\pm \tau ^{2}\alpha $ of an
integer $\alpha$.
Obviously  $|\mu \alpha | ^{2s}=|\alpha | ^{2s}\ $
(for any unit $\mu $ of $K$). Since
$$ \overline{-\alpha }=-\overline{\alpha },\quad \overline{-\tau \alpha }=
-\overline{\tau \alpha },\quad \overline{-\tau ^{2}\alpha }=
-\overline{\tau ^{2}\alpha }, $$
so
$$ \frac{\overline{\alpha }}{|\alpha | ^{2s}}+\frac{\overline
{-\alpha }}{|-\alpha | ^{2s}}+ \frac{\overline{\tau \alpha }}{|\tau \alpha |
^{2s}}+\frac{\overline{-\tau \alpha }}{|-\tau \alpha | ^{2s}}+\frac{\overline
{\tau  ^{2}\alpha }}{|\tau  ^{2}\alpha | ^{2s}}+\frac{\overline
{-\tau  ^{2}\alpha }}{|-\tau  ^{2}\alpha | ^{2s}}=0  .$$
 Therefore we have
$$\quad  \sum _{\alpha \in L-\{0\}} \frac{\overline{\alpha }}{|\alpha | ^{2s}}
=0  . $$
Similarly, by
$$ \frac{\overline{-\alpha }}{-\alpha }=\frac{\overline{\alpha }}{\alpha } \qquad  \frac{\overline{\tau \alpha  }}{\tau \alpha }=\tau \cdot
\frac{\overline{\alpha }}{\alpha } \qquad  \frac{\overline
{-\tau \alpha }}{-\tau \alpha }=\tau \cdot \frac{\overline{\alpha }}
{\alpha }, $$
$$\quad \frac{\overline{\tau  ^{2}\alpha  }}{\tau  ^{2}\alpha }=\tau  ^{2}\cdot
\frac{\overline{\alpha }}{\alpha } \qquad  \frac{\overline{-\tau  ^{2}\alpha  }}
{-\tau  ^{2}\alpha }=\tau  ^{2}\cdot
\frac{\overline{\alpha }}{\alpha } \ ,  $$
and $\quad 1+\tau +\tau ^{2}=0,\ $ so we know
$$ \frac{\overline{-\alpha }}{-\alpha }+\frac{\overline{\alpha }}{\alpha }+
\frac{\overline
{\tau \alpha  }}{\tau \alpha }+\frac{\overline{-\tau \alpha }}{-\tau \alpha
}+ \frac{\overline{\tau  ^{2}\alpha  }}{\tau  ^{2}\alpha }+\frac{\overline{-\tau  ^{2}\alpha  }}
{-\tau  ^{2}\alpha }=2(1+\tau +\tau ^{2})\cdot \frac{\overline{\alpha }}{\alpha
}=0, $$
 Thus we have
$$\quad \sum _{\alpha \in L-\{0\}} \frac{1}{|\alpha | ^{2s}}
\cdot  \frac{\overline{\alpha }} {\alpha }=0 .$$
Therefore, when $Re(s)$ is sufficiently large,
by the series $(3.4 ),\ $
defined above we have
\begin{align*}  \sum _{\alpha \in L} \frac{\overline{z}+\overline{\alpha }}
{|z+\alpha | ^{2s}}&=\psi (z, s, L)+\sum _{\alpha \in L-\{0\}} \frac{\overline
{\alpha }}{|\alpha | ^{2s}}-sz \sum _{\alpha \in L-\{0\}} \frac{1}
{|\alpha | ^{2s}}\cdot \frac{\overline{\alpha }}{\alpha }+(1-s)
\sum _{\alpha \in L-\{0\}} \frac{\overline{z}}{|\alpha | ^{2s}}\\
&=\psi (z, s, L)+(1-s)\overline{z}\sum _{\alpha \in L-\{0\}}\frac{1}
{(N(\alpha )) ^{s}}.
\end{align*}
That is
$$ \sum _{\alpha \in L} \frac{\overline{z}+\overline{\alpha }}
{|z+\alpha | ^{2s}}=\psi (z, s, L)+6\overline{z}(1-s)\zeta _{K}(s) \qquad  (3.6)$$
where $\zeta _{K}(s)$is the Dedekind  Zeta-function of the
number field $K=\Q(\sqrt{-3})$
Since $\psi (z, s, L)$ is  an  analytic function in the
area $\quad Re(s)> \frac{1}{2}$, and $\zeta _{K}(s)$
is  an  analytic function for $\quad Re(s)>1$,
so the right side of   formula $(3.6)$ gives an
analytic extension for the series
 $\quad \sum _{\alpha \in L} \frac{\overline{z}+
 \overline{\alpha }}{|z+\alpha | ^{2s}}\quad $
(to the area $\quad Re(s) >1\quad $). Now transform the
right side of  formula $(3.3 )$ as follows:
\begin{align*}
\sum _{v} \frac{3\overline{v}\overline{\Delta  }+
3\overline{\beta }+\overline{\Delta  }}{N(3v\Delta  +3\beta +\Delta  ) ^{s}}
&=\sum _{v} \frac{3\overline{\Delta  } \left(\overline{v}+\frac{3\overline{\beta }+
\overline{\Delta  }}{3\overline{\Delta  }} \right)}{(N(3\Delta  )) ^{s}
\left(N \left(v+\frac{3\beta +\Delta  }{3\Delta  } \right) \right) ^{s}}\\
&=\frac{3\overline{\Delta  }}{N(3\Delta  ) ^{s}} \sum _{v} \frac{\overline{v}+
\frac{3\overline{\beta }+\overline{\Delta  }}{3\overline{\Delta  }}}{N \left(v+
\frac{3\beta +\Delta  }{3\Delta  } \right) ^{s}} \quad .
\end{align*}
By the analytic extension of $(3.6)$  we could
 obtain the analytic extension of
 $$\sum \limits _{v} \frac{\overline{v}+
\frac{3\overline{\beta }+\overline{\Delta  }}{3\overline{\Delta  }}}{N \left(v+
\frac{3\beta +\Delta  }{3\Delta  } \right) ^{s}},  $$
and hence  obtain  the analytic extension of $L_{D}(s)$
as the following:
\begin{align*}
L_{D}(s)&=\sum _{\beta } \left(\pm \left(\frac{3\beta }{D} \right) _{6} \right)
\sum _{v} \frac{3\overline{v}\overline{\Delta  }+3\overline{\beta }+
\overline{\Delta  }}{N(3v\Delta  +3\beta +\Delta  )^{s}}\\
&=\sum _{\beta }\left(\pm \left(\frac{3\beta }{D} \right) _{6} \right)\cdot
\frac{3\overline{\Delta  }}{N(3\Delta  ) ^{s}}\cdot \sum _{v} \frac{\overline{v}+
\frac{3\overline{\beta }+\overline{\Delta  }}{3\overline{\Delta  }}}{N \left(v+
\frac{3\beta +\Delta  }{3\Delta  } \right) ^{s}}\\
&=\frac{3\overline{\Delta  }}
{N(3\Delta  ) ^{s}}\cdot \sum _{\beta } \left(\pm \left(\frac{3\beta }
{D} \right) _{6} \right) \left[\psi (\frac{\beta }{\Delta
}+\frac{1}{3}, s, L)+6(\frac{\overline{\beta }}{\overline{\Delta
}}+\frac{1}{3})(1-s)\zeta _{K}(s) \right]  .
\end{align*}
By class number formula of imaginary quadratic field
( see e.g. [Wa], [Zhang] ) we  obtain
$$ \lim _{s\rightarrow 1}(s-1)\zeta _{K}(s)=\frac{2\pi }{6\sqrt{3}}\cdot
h(K)=\frac{2\pi }{6\sqrt{3}}=\frac{\pi }{3\sqrt{3}} . $$
Thus let  $s\rightarrow 1,\ $   we obtain
$$ L_{D}(1)=\frac{3\overline{\Delta  }}{N(3\Delta  )}\sum _{\beta }
\left(\pm \left(\frac{3\beta }{D} \right) _{6} \right)\left[\psi (\frac{\beta
}{\Delta  }+\frac{1}{3},1,L)-\frac{2\pi }{\sqrt{3}}(\frac{\overline{\beta }}
{\overline{\Delta  }}+\frac{1}{3}) \right]   . $$
Also we have $\psi (z,1,L)=\zeta (z,L),\ $   therefore
$$ L_{D}(1)=\frac{1}{3\Delta  } \sum _{\beta }
\left(\pm \left(\frac{3\beta }{D} \right) _{6} \right) \left[\zeta (\frac{\beta }
{\Delta  }+\frac{1}{3},L)-\frac{2\pi }{\sqrt{3}}(\frac{\overline{\beta }}
{\overline{\Delta  }}+\frac{1}{3}) \right]  \qquad (3.7)$$
By results in [ Ste ]  we know , the Weiestrass $ \wp -$function
corresponding to the period lattice $L=O_{K}=\Z+\Z _{\tau }$
  is  $\Omega ^{2} \wp (\Omega z),\ $ where function $\wp (z)$
satisfies  $\wp '(z) ^{2}=4\wp (z) ^{3}-1,\ $  and
the corresponding period lattice  is  $\Omega O_{K}\ $
($\Omega =3.059908 \cdots  $). And  we have  the following
  formula :
$$ \zeta (\alpha _{1}+\alpha _{2})=\zeta (\alpha _{1})+\zeta (\alpha _{2})+
\frac{\Omega }{2}\cdot \frac{\wp '(\Omega \alpha _{1})-
\wp '(\Omega \alpha _{2})}{\wp (\Omega \alpha _{1})-\wp
(\Omega \alpha _{2})} . \qquad  (3.8)$$
 Thus
$$ L_{D}(1)=\frac{1}{3\Delta  } \sum _{\beta }\left(\pm
\left(\frac{3\beta }{\Delta  } \right) _{6} \right)\cdot $$
$$\qquad \left[\zeta (\frac{\beta } {\Delta  })+\zeta
(\frac{1}{3})+
\frac{\Omega }{2}\cdot \frac{\wp '(\frac{\beta \Omega }
{\Delta  })-\wp '(\frac{\Omega }{3})}{\wp (\frac{\beta \Omega }
{\Delta  })-\wp (\frac{\Omega }{3})}-\frac{2\pi }{\sqrt{3}}
(\frac{\overline{\beta }} {\overline{\Delta  }}+\frac{1}{3})
\right]   . $$
Also by  [ Ste ]  we know  $\quad \wp '(\frac{\Omega }{3})=
-\sqrt{3},\ $  $\wp (\frac{\Omega }{3})=1,\ $
so  we obtain  the following Proposition :

{\bf  Proposition  3.1 }.\quad  Let   $D\equiv 1(mod\ 3)$
be a primitive integer  of $ O_{K}.\ $ Then for the
elliptic curve $E_{D}:\ y^{2}=x^{3}+2^{4}D,\ $ the $L-$series
$L_{D}(s)=\sum \limits _{\sigma \equiv 1
(mod\ 3)}\left(\frac{D}{\sigma } \right) _{6} \cdot \frac{\overline{\sigma }}
{(N(\sigma )) ^{s}}$ could be analytic extended via
the series $\psi (z, s, L)$ and Dedekind  Zeta-function
$\zeta _{K}(s)$ , and we have
$$L_{D}(1)=\frac{1}{3\Delta  } \sum _{\beta }
\left(\pm \left(\frac{3\beta }{D} \right) _{6} \right)\cdot $$
$$\left[\zeta (\frac{\beta } {\Delta  })+
\zeta (\frac{1}{3})+\frac{\Omega }{2}\cdot \frac{\wp '(\frac{\beta \Omega }
{\Delta  })+\sqrt{3}}{\wp (\frac{\beta \Omega }{\Delta  })-
1}-\frac{2\pi }{\sqrt{3}}(\frac{\overline{\beta }}
{\overline{\Delta  }}+\frac{1}{3}) \right]  \ (3.9)$$

{\bf   Remark } 3.2 .\quad By formula (3.9) in the
 above Proposition  3.1 ,  we could in particular obtain
 the corresponding result in [Ste] on $L-$series $L_{D}(s)$
 for elliptic curves
$E_{D}:\ y^{2}=x^{3}-2^{4}3^{3}D^{2}$
( with $D$  rational  integer ).

    \par \vskip 0.4cm
Now we turn to prove Theorem 1.6 . Under the assumption of
Theorem1.6  we have the following lemma by the definition of
the $L-$series:
    \par \vskip 0.4cm

{\bf  Lemma   3.1 }. $\quad L_{S}(\overline{\psi } _{D_{T} ^{2}}, s)=\left\{
\begin{array}{l}  L(\overline{\psi } _{D_{T} ^{2}}, s)  \ ,  \qquad \hbox{if }\quad  \prod \limits _{\pi _{k} \in S} \pi _{k}=D_{T} \quad ;\\
L(\overline{\psi } _{D_{T} ^{2}}, s) \prod \limits _{\pi _{k}|\widehat{D}_{T}}
\left(1-\left(\frac{D_{T}}{\pi _{k}} \right) _{3}\cdot \frac{\overline{\pi _{k}}}
{(\pi _{k}\overline{\pi _{k}}) ^{s}} \right)  \qquad  \hbox{ otherwise  }.
\end{array}
\right.$
  \par \vskip 0.4cm

{\bf  Proof of Theorem 1.6 . }.\quad
For the elliptic curve
$\quad E_{D_{T} ^{2}}:\ y^{2}=x^{3}-2^{4}3^{3}D_{T} ^{2}\ $,
assume its period lattice  is  $L_{T}=\omega _{T}O_{K} \ .$
Since the class number of  $K=\Q (\sqrt {-3})$ is a $h_{K}=1,\ $
by  [Sil\ 2]  we know  that all elliptic curves defined
over  $\C$ with complex multiplication ring $O_{K}$
are $\C-$isomorphic each other. So their  period lattices
are Homothetic each other. We know the elliptic curve
corresponding to the lattice $O_{K}$, denoted by $\C/O_{K}$,
has complex multiplication ring $O_{K}.\ $
Therefore every elliptic curve $E$ over $\C$ with
complex multiplication ring $O_{K}$ has  period lattice $L$
Homothetic to  $O_{K}$, i.e.  we always have $L=\beta O_{K}$
(for some  $\beta \in \C ^{\times}$).
Thus for the above elliptic curve $E_{D_{T}^{2}}$
and its period lattice $L_{T}=\omega _{T} O_{K},\ $
we may assume  $\omega _{T}=\alpha _{T}\omega ,\ $
$\alpha _{T}\in \C^{\times }$ (In fact, it's easy to see that
$\omega _{T}={\omega } ({2\sqrt{3}\sqrt[3] {D_{T}}})^{-1},\ $
  i.e.
$\alpha _{T}={1}({2\sqrt{3}\sqrt[3] {D_{T}}})^{-1}$).
By  [Ste, P.125]  we know that the conductor of
 $\psi _{D_{T}^{2}}$ is $\sqrt{-3}D_{T}$ or $3D_{T} \ .$
 Therefore, in Proposition $(A)$ above, putting
  $k=1,\ $ $\rho =\frac{\omega _{T}}{(3D)},\ $ $\h=O_{K},\ $
 $\g =(3D),\ $ $\phi =\psi _{D_{T} ^{2}},\ $   we obtain
$$ \frac{\overline{\rho }}{|\rho | ^{2s}} L_{\g}(\overline{\psi } _{D_{T}
^{2}}, s)=\sum _{\flat \in \mathbf{B}} H_{1} \left(
\frac{\psi _{D_{T} ^{2}}(\flat )\omega _{T}}{3D}, 0, s, L_{T} \right), \quad
(Re(s)> {3}/{2}) $$
Since the conductor $f$ of $\psi _{D_{T} ^{2}}$
divides $\ (3D)=\g,\ $ so by Lemma B above  we know  the
ray class of $K$ to the cycle (modulo)  $\g$ is
 $K((E_{D_{T} ^{2}}) _{(3D)}),\ $  In particular, we have
$$ \left(O_{K}/(3D) \right) ^{\times }/\mu _{6} \cong Gal \left(K((E_{D_{T} ^{2}})
_{(3D)})/K \right)  \qquad  (via \quad Artin \hbox{ map }) $$
where $\mu _{6}=\{\pm 1,\pm \tau ,\pm \tau ^{2} \}$ and
$\mu _{6}\cong (O_{K}/(3)) ^{\times } \ .$ So we may take  the
set $\mathbf{B}=\{(3c+D):\ c \in \mathcal{C}\},\ $ where
$\mathcal{C}$ is a reduced residue system of $O_{K}$ modulo $D$
 i.e. a representative system for
  $\mathcal{C}$is $(O_{K}/(D)) ^{\times }$. Thus
$$ \frac{\overline{\rho }}{|\rho | ^{2s}} L_{\g}(\overline{\psi } _{D_{T}
^{2}}, s)=\sum _{c \in \mathcal{C}} H_{1} \left(
\frac{\psi _{D_{T} ^{2}}(3c+D)\omega _{T}}{3D},0, s,\omega _{T}O_{K}
\right) ,\quad (Re(s)> {3}/{2})$$
Note that $\quad H_{1}(z,0,1,L) \ $ could be analytically
continued by the Eisenstein $\ E^{*}-$function
( see [ We\ 2 ] ):
$$ H_{1}(z,0,1,L)=E_{0,1} ^{*}(z,L)=E_{1}^{*}(z,L)  . $$
Hence  we have
$$ \frac{1}{\rho }L_{\g}(\overline{\psi }_{D_{T} ^{2}},1)=
\sum _{c \in \mathcal{C}} E_{1} ^{*} \left(\frac{\psi _{D_{T} ^{2}}
(3c+D)\omega _{T}}{3D},\omega _{T}O_{K} \right),\quad (3.10)$$
 that is
$$ \frac{3D}{\alpha _{T}\ \omega }L_{\g}(\overline{\psi }_{D_{T} ^{2}},1)=
\sum _{c \in \mathcal{C}} E_{1} ^{*} \left(\frac{\psi _{D_{T} ^{2}}
(3c+D)\alpha _{T}\ \omega }{3D},\alpha _{T}\ \omega O_{K} \right)  . $$
Since $D\equiv 1(mod\ 6),\ $ so $3c+D\equiv 1(mod\ 3)\ $
 for any $c \in \mathcal{C}.\ $ Thus by the definition of
 $\psi _{D_{T} ^{2}}$  and the cubic reciprocity, we have
\begin{align*}
\psi _{D_{T} ^{2}}(3c+D)&=\overline{\left(\frac{D_{T}}{3c+D} \right) _{3}}(3c+D)
=\overline{\left(\frac{3c+D}{D_{T}} \right) _{3}}(3c+D)\\
&=\overline{\left(\frac{3c}{D_{T}} \right) _{3}}(3c+D)=
\left(\frac{3c}{D_{T}} \right) _{3} ^{2}(3c+D)
\end{align*}
 Note that $\quad L_{\g}(\overline{\psi }_{D_{T} ^{2}},1)=
L_{S}(\overline{\psi }_{D_{T} ^{2}},1) ,\ $ so by
(3.10) we have
$$ \frac{3D}{\alpha _{T}\ \omega } L_{S}(\overline{\psi }_{D_{T}
^{2}},1)=\sum _{c \in \mathcal{C}}E_{1} ^{*}\left(\left(\frac{c\ \omega }{D}+
\frac{\omega }{3}\right)\alpha _{T} \left(\frac{3c}{D_{T}} \right) _{3} ^{2},
\quad \alpha _{T} \ \omega O_{K} \right) \ (3.11)$$
Let $\lambda =\alpha _{T} \left(\frac{3c}{D_{T}}
\right) _{3} ^{2},\ $ by   formula
$E_{1} ^{*}(\lambda z,\lambda L)=
\lambda ^{-1}E_{1} ^{*}(z,L),\ $  we obtain
\begin{align*}
E_{1} ^{*}\left(\left(\frac{c\ \omega }{D}+\frac{\omega }{3}
\right) \alpha _{T} \left(\frac{3c}{D_{T}} \right) _{3} ^{2}\ ,
\ \alpha _{T} \omega O_{K} \right)&= E_{1} ^{*}\left(\left(
\frac{c\ \omega }{D}+\frac{\omega }{3} \right)\lambda \ ,
\ \lambda L_{\omega }\right)\\
\hskip -2cm & =\lambda ^{-1} E_{1} ^{*}\left(\left(
\frac{c\ \omega }{D}+\frac{\omega }{3}\right)\ ,
\ L_{\omega } \right) \\
\hskip -2cm &=\alpha _{T} ^{-1} \left(
\frac{3c}{D_{T}} \right) _{3} E_{1} ^{*}
\left(\left(\frac{c\ \omega }{D}+
\frac{\omega }{3} \right) \ ,\ L_{\omega }\right).
\end{align*}
So by $(3.11)$, we have
 \begin{align*}
\frac{3D}{\omega } L_{S} (\overline{\psi } _{D_{T} ^{2}},1)&=
\sum _{c\ \in \mathcal{C}} \left(\frac{3c}{D_{T}} \right) _{3}
E_{1} ^{*} \left(\left(\frac{c\ \omega }{D}+
\frac{\omega }{3} \right) \ ,
\ L_{\omega }\right)\\
\hskip -2cm &=\left(\frac{3}{D_{T}} \right) _{3}\
\sum _{c\ \in \mathcal{C}}
\left(\frac{c}{D_{T}} \right) _{3}
E_{1} ^{*} \left(\left(\frac{c\ \omega }{D}+
\frac{\omega }{3} \right) \ , \ L_{\omega }\right) ,
\end{align*}
$$\frac{D}{\omega } \left(\frac{9}{D_{T}} \right) _{3}\
L_{S} (\overline{\psi } _{D_{T} ^{2}},1)=\frac{1}{3}\sum _{c\ \in \mathcal{C}}
\left(\frac{c}{D_{T}} \right) _{3}E_{1} ^{*} \left(\left(\frac{c\ \omega }{D}+
\frac{\omega }{3} \right) \ ,\ L_{\omega }\right), \quad (3.12)$$
By  [Go-Sch, Prop.1.5] we know
$$ E_{1} ^{*}(z,L)=\zeta (z,L)-zs_{2}(L)-\overline{z}A(L) ^{-1} $$
where
$$ \zeta (z,L)=\frac{1}{z}+\sum _{\alpha \in L-\{0\}}
\left(\frac{1}{z-\alpha }
+\frac{1}{\alpha }+\frac{z}{\alpha ^{2}} \right) $$
is  Weierstrass $  \ Zeta-$function, an odd function. For
$$S_{2}(L)=\lim \limits _{s\rightarrow 0 \ s> 0}\
\sum \limits _{\alpha \in L-\{0\}}
\alpha ^{-2} |\alpha | ^{-2s},\ $$  we have
$$\eta (\alpha ,L)=\zeta (z+\alpha ,L)-\zeta (z,L)\ ,\ \eta(\alpha ,L)=
\alpha S_{2}(L)+\overline{\alpha }A(L)^{-1}  \ ,  $$
( for any $\alpha \in L$),  $\eta $ is a quasi-period
map corresponding to $L$. The Weierstrass
$\ \wp -$function corresponding to the
  period lattice $L_{\omega }=\omega \ O_{K}\ $ is
$$\wp (z,L_{\omega }):\ \wp '(z)^{2}=4\wp (z)^{3}-1, $$
and in this case
$$ A(L_{\omega }) =\frac{\overline{\omega }(\omega \tau )-
\omega (\overline{\omega \tau })}{2\pi I}
=\frac{\omega ^{2}(\tau  -\overline{\tau })}{2\pi I}
=\frac{\sqrt{3} \omega ^{2}}{2\pi }  . $$
Also obviously we have
$\frac{\omega }{2} \not\in L_{\omega },\ $  therefore
by  [Sil\ 2 p.41] we know
$$ \eta (\omega ,L_{\omega })=2\zeta
(\frac{\omega }{2},L_{\omega })  \qquad
 \eta (\omega , L_{\omega })=\omega S_{2}(L_{\omega })+
\overline{\omega }A(L_{\omega }) ^{-1} . $$
Thus
$$ 2\zeta (\frac{\omega }{2},L_{\omega })=\omega S_{2}(L_{\omega })+
\overline{\omega } \cdot \frac{2\pi }{\sqrt{3}\omega ^{2}}=
\omega S_{2}(L_{\omega })+\frac{2\pi }{\sqrt{3}\omega }. $$
So
$$ S_{2}(L_{\omega })=\frac{2}{\omega }\zeta (\frac{\omega }{2},
L_{\omega })-\frac{2\pi }{\sqrt{3}\omega ^{2}}  . $$
  Therefore
$$ E_{1} ^{*}(z,L_{\omega })=\zeta (z,L_{\omega })-\frac{2z}{\omega }
\zeta (\frac{\omega }{2}, L_{\omega })+\frac{2\pi z}{\sqrt{3}\omega ^{2}}
-\frac{2\pi }{\sqrt{3}\omega ^{2}} \overline{z}  . $$
Put $z=\frac{c \omega }{D}+\frac{\omega }{3},\ $  we obtain
\begin{align*}
E_{1} ^{*}(\frac{c \omega }{D}+\frac{\omega }{3},L_{\omega })&=
\zeta (\frac{c \omega }{D}+\frac{\omega }{3},L_{\omega })-
\frac{2}{\omega } \left(\frac{c \omega }{D}+\frac{\omega }{3} \right)\zeta
(\frac{\omega }{2},L_{\omega })\\
&+\frac{2\pi }{\sqrt{3}\omega ^{2}}
\left(\left(\frac{c \omega }{D}+\frac{\omega }{3} \right)-
\overline{\left(\frac{c\ \omega }{D}+\frac{\omega }{3} \right)} \right)\\
&=\zeta (\frac{c \omega }{D}+\frac{\omega }{3},L_{\omega })-
2(\frac{c}{D}+\frac{1}{3})\zeta (\frac{\omega }{2},L_{\omega })+
\frac{2\pi }{\sqrt{3}\omega }(\frac{c}{D}-
\frac{\overline{c}}{\overline{D}}).
\end{align*}
  That is
\begin{align*}
E_{1} ^{*}(\frac{c \omega }{D}+\frac{\omega }{3},L_{\omega })&=
\zeta (\frac{c \omega }{D}+\frac{\omega }{3},L_{\omega })\\
&-2(\frac{c}{D}+\frac{1}{3})\zeta (\frac{\omega }{2},L_{\omega })+
\frac{2\pi }{\sqrt{3}\omega }(\frac{c}{D}-\frac{\overline{c}}{\overline{D}})
\quad  \qquad  (3.13)
\end{align*}
Now let us show $\zeta (\frac{\omega }{3},L_{\omega })$ and
 $\zeta (\frac{\omega }{2}, L_{\omega })$ are equal.
 In fact, by formula (3.2) in [Ste, P.126],
for  any  rational  integer $D(3\dag D),\ $ we know
\begin{align*}
L_{D}(s)&=\sum _{\sigma } \left(\frac{\sigma }{D} \right) _{3}
\sum _{u \in \Q(\tau )}
\frac{3\Delta  \overline{u }+3\overline{\sigma }+\Delta  }{N(3\Delta  u
+3\sigma +\Delta  ) ^{s}}\\
&=\frac{3\Delta  }{N(3\Delta  ) ^{s}} \left\{\sum _{\sigma
} \left(\frac{\sigma }{D} \right) _{3} \psi (\frac{\sigma }{\Delta  }+
\frac{1}{3}, s)+6(1-s)\zeta _{K}(s)\sum _{\sigma
} \left(\frac{\sigma }{D} \right) _{3} \left(\frac{\overline{\sigma }}{\Delta  }
+\frac{1}{3} \right) \right\} \ .
\end{align*}
Let $s\rightarrow 1,\ $ since
$\quad \lim \limits _{s\rightarrow 1}\  (s-1)
\zeta _{K}(s)=\frac{\pi }{3\sqrt{3}},\ $ and
 $\psi (z, s, L)$ is convergent to $\xi (z,L)$ uniformly
( when $s\rightarrow 1$), so
$$ L_{D}(1)=\frac{1}{3\Delta  } \left[\sum _{\sigma } \left(\frac{\sigma }{D}
\right) _{3} \xi (\frac{\sigma }{\Delta  }+\frac{1}{3})-\frac{2\pi }{\sqrt{3}}
\sum _{\sigma } \left(\frac{\sigma }{D}
\right) _{3}  (\frac{\overline{\sigma }}{\Delta  }+
\frac{1}{3}) \right] $$
Put $\quad u ={\rho }/{3\Delta  },\ $ $\rho '=\rho +3\Delta
,\ $   $\quad \rho =3\sigma +\Delta  $, then we have
\begin{align*}
L_{D}(1)&=\frac{1}{3\Delta  } \left[\sum _{\sigma } \left(\frac{\sigma }{D}
\right) _{3} \xi (\frac{\rho }{3 \Delta  })-\frac{2\pi }{\sqrt{3}}
\sum _{\sigma } \left(\frac{\sigma }{D}
\right) _{3}\left(\frac{\overline{\rho }}{ 3 \Delta  }\right) \right]  \\
&=\frac{1}{3\Delta  } \left[\sum _{\sigma '} \left(\frac{\sigma '}{D}
\right) _{3} \xi (\frac{\rho '}{3 \Delta  })-\frac{2\pi }{\sqrt{3}}
\sum _{\sigma '} \left(\frac{\sigma '}{D}
\right) _{3}\left(\frac{\overline{\rho '}}{ 3 \Delta  }\right) \right] \\
&=\frac{1}{3\Delta  } \left[\sum _{\sigma } \left(\frac{\sigma }{D}
\right) _{3} \xi (u +1)-\frac{2\pi }{\sqrt{3}}
\sum _{\sigma } \left(\frac{\sigma }{D}
\right) _{3}\left(\frac{\overline{\rho }}{ 3 \Delta  }+1\right) \right] \\
&=\frac{1}{3\Delta  } \left[\sum _{\sigma } \left(\frac{\sigma }{D}
\right) _{3} \xi (u+1)-\frac{2\pi }{\sqrt{3}}
\sum _{\sigma } \left(\frac{\sigma }{D}
\right) _{3}\left(\frac{\overline{\rho }}{ 3 \Delta  }\right)-\frac{2\pi }
{\sqrt{3}}\sum _{\sigma } \left(\frac{\sigma }{D}\right) _{3} \right]
\end{align*}
 therefore   we obtain
$$ \frac{1}{3\Delta  } (\xi (u +1)-\xi (u ))= \frac{1}{3\Delta  }\cdot
\frac{2\pi }{\sqrt{3}} \, $$
 that is $\quad \xi (u +1)=\xi (u )+\frac{2\pi }{\sqrt{3}}.\ $
Similarly we  could obtain
$\quad \xi (u +\tau )=\xi (u )+\frac{2\pi }{\sqrt{3}}
\overline{\tau } \quad .$   Therefore  we have
$$\quad \xi (u +1)=\xi (u )+\frac{2\pi }{\sqrt{3}}  \ ,  \qquad \xi (u +\tau )=\xi (u )+\frac{2\pi }{\sqrt{3}} \overline{\tau }
 \qquad (3.14) $$
where $\xi (z)$ is the  Weierstrass
$\ Zeta-$function with period lattice $\ L=\Z+\Z \tau =O_{K}$.
    \par \vskip 0.3cm

In formula $(3.14)$, Let $ \quad u=-{1}/{2}\ $ or
 $-\frac{\tau }{2}.\ $ Since
$\xi (u)$ is an odd function , we obtain
$$ \xi (\frac{1}{2})=\frac{\pi }{\sqrt{3}}  \qquad    \xi (\frac{\tau }{2})
=\frac{\pi }{\sqrt{3}}\overline{\tau }  \ ,  $$
 therefore
$$ \zeta (\frac{\omega }{2},L_{\omega })=\frac{1}{\omega }\xi (\frac{1}{2},
O_{K})=\frac{\pi }{\sqrt{3}\ \omega }  \ ,  $$
 that is
$$ \qquad  \zeta (\frac{\omega }{2},L_{\omega })=\frac{\pi
}{\sqrt{3}\ \omega }  \qquad   \qquad  (3.15) $$
Also  by  [Ste, P.127] we know,
$$\wp (\frac{\omega }{3},L_{\omega })=1  \ ,  \quad \wp '(\frac{\omega
}{3},\ L_{\omega })=-\sqrt{3} \quad  \qquad  (3.16) $$
so by the  formulae
$$ \wp ''(z)=6\wp (z) ^{2} -\frac{1}{2}g_{2}  \qquad  2\zeta (2z)-
4\zeta (z)=\frac{\wp ''(z)}{\wp '(z)}  \ ,  $$
( see [Law], P.182),
let $\quad z={\omega }/{3},\ $  we obtain
$\quad\wp ^{\prime \prime}({\omega }/{3}, L_{\omega })=
6\wp ({\omega }/{3}) ^{2}=6, $
 so
$$ 2\zeta ({2\omega }/{3}, L_{\omega })
-4\zeta ({\omega }/{3}, L_{\omega })
= -{6}/{\sqrt{3}}=-2\sqrt{3}, $$
 i.e. $ \  \zeta ({2\omega }/{3}, L_{\omega })
 -2\zeta ({\omega }/{3}, L_{\omega })=-\sqrt{3} \ .$
  On the other hand, in  formula $(3.14)$, let
  $\ u=-{1}/{3},\ $  we obtain
  $\ \xi (-{1}/{3}+1)= \xi (-{1}/{3})+{2\pi }/{\sqrt{3}},\ $
 i.e.  $\ \xi ({2}/{3})+\xi ({1}/{3})={2\pi }/{\sqrt 3}  .$
 Also we have
$$ \zeta ({2\omega }/{3}, L_{\omega })=
{\omega }^{}-1 \xi ({2}/{3}),
 \qquad   \zeta ({\omega }/{3}, L_{\omega })=
 {\omega }^{-1}\xi ({1}/{3}), $$
so
$$ \omega \zeta (\frac{2\omega }{3},L_{\omega })+\omega \zeta (\frac{\omega }
{3},L_{\omega })=\frac{2\pi }{\sqrt{3}}  \ ,  $$
which gives
$$ \left\{
\begin{array}{l}  \zeta ({2\omega }/{3}, L_{\omega })+
\zeta ({\omega }/{3}, L_{\omega })=
{2\pi }/({\sqrt{3} \omega }) , \\
\zeta ({2\omega }/{3}, L_{\omega })-2\zeta
({\omega }/{3}, L_{\omega })=-\sqrt{3}.
\end{array}
\right. $$
This gives the solution
$$ \zeta (\frac{\omega }{3},L_{\omega })=\frac{2\pi }{3\sqrt{3} \omega
}+\frac{1}{\sqrt{3}}  \ ,  \quad \zeta (\frac{2\omega }{3},L_{\omega
})=\frac{4\pi }{3\sqrt{3} \omega }-\frac{1}{\sqrt{3}} \quad .\quad
(3.17)$$
 Also  by  the formula ( see [Law] ):
$$ \zeta (z_{1}+z_{2},L_{\omega })=\zeta (z_{1},L_{\omega })+
\zeta (z_{2},L_{\omega })+\frac{1}{2} \frac{\wp '(z_{1})-\wp '(z_{2})}
{\wp (z_{1})-\wp (z_{2})}  \ ,  $$
 we obtain
\begin{align*}
\zeta (\frac{c\ \omega }{D}+\frac{\omega }{3},L_{\omega })&=
\zeta (\frac{c\ \omega }{D},L_{\omega })+\zeta (\frac{\omega }{3},L_{\omega })
+\frac{1}{2} \frac{\wp '(\frac{c\ \omega }{D},L_{\omega })-
\wp '(\frac{\omega }{3},L_{\omega })}{\wp (\frac{c\ \omega }{D},L_{\omega })-
\wp (\frac{\omega }{3},L_{\omega })}\\
& = \zeta (\frac{c\ \omega }{D},L_{\omega })
+ \frac{2\pi }{3\sqrt{3} \omega }+\frac{1}{\sqrt{3}}+
\frac{1}{2} \frac{\wp '(\frac{c\ \omega }{D},L_{\omega })+\sqrt{3}}
{\wp (\frac{c\ \omega }{D},L_{\omega }) -1}
\end{align*}
substitute this into $(3.13)$, we obtain
\begin{align*}
E_{1} ^{*}(\frac{c \omega }{D}+\frac{\omega }{3},L_{\omega })&=
\zeta (\frac{c\ \omega }{D},L_{\omega })
+ \frac{2\pi }{3\sqrt{3} \omega }+\frac{1}{\sqrt{3}}+
\frac{1}{2} \frac{\wp '(\frac{c\ \omega }{D},L_{\omega })+\sqrt{3}}
{\wp (\frac{c\ \omega }{D},L_{\omega }) -1} \\
&-2(\frac{c}{D}+\frac{1}{3})\frac{\pi }{\sqrt{3} \omega }+
\frac{2\pi }{\sqrt{3}\omega
}(\frac{c}{D}-\frac{\overline{c}}{\overline{D}})\\
&=\zeta (\frac{c\ \omega }{D},L_{\omega })+\frac{1}{2}
\frac{\wp '(\frac{c\ \omega }{D},L_{\omega })+\sqrt{3}}
{\wp (\frac{c\ \omega }{D},L_{\omega }) -1}+\frac{1}{\sqrt{3}}-
\frac{2\pi }{\sqrt{3}\ \omega }\cdot \frac{\overline{c}}{\overline{D}}
\end{align*}
Now substitute this into $(3.12)$, we have
\begin{align*}
&\frac{D}{\omega } \left(\frac{9}{D_{T}} \right) _{3}\
L_{S} (\overline{\psi } _{D_{T} ^{2}},1)\\
&=\frac{1}{3}\sum _{c\ \in \mathcal{C}}
\left(\frac{c}{D_{T}} \right) _{3}
\left[\zeta (\frac{c\ \omega }{D},L_{\omega })
+\frac{1}{2} \frac{\wp '(\frac{c\ \omega }{D},L_{\omega })+\sqrt{3}}
{\wp (\frac{c\ \omega }{D},L_{\omega }) -1}+\frac{1}{\sqrt{3}}-
\frac{2\pi }{\sqrt{3}\ \omega }\cdot \frac{\overline{c}}
{\overline{D}} \right] .
\end{align*}
Since  $\quad D=\pi _{1} \cdots \pi _{n}\ $ with
 $\pi _{k}\equiv 1 (mod\ 6),\ $
 so we may choose the representatives $\mathcal{C}\ $
for $\quad \left(O_{K}/(D)\right) ^{\times } \ $
in such a way that  $\ -c \in \mathcal{C} \ $ when
$c \in \mathcal{C}$.
Obviously  $\ \left({-c}/{D_{T}} \right) _{3} =
\left({c}/{D_{T}} \right) _{3} . $
Also since  $\zeta (z,L_{\omega })$ and
$\wp '(z,L_{\omega })$ are odd functions , and
$\wp (z,L_{\omega })$ is even function, so
\begin{align*}
\sum _{c\ \in \mathcal{C}}
\left(\frac{c}{D_{T}} \right) _{3} \zeta (\frac{c\ \omega }{D},L_{\omega })
&=\sum _{c\ \in \mathcal{C}} \left(\frac{c}{D_{T}} \right) _{3}
\frac{\wp '(\frac{c\ \omega }{D},L_{\omega })}
{\wp (\frac{c\ \omega }{D},L_{\omega }) -1}\\
&=\sum _{c\ \in \mathcal{C}} \left(\frac{c}{D_{T}} \right) _{3}\cdot
\frac{\overline{c}}{\overline{D}}=0
\end{align*}
  Therefore  , $$ \frac{D}{\omega }\left(\frac{9}{D_{T}} \right) _{3} L_{S}(\overline
{\psi } _{D_{T} ^{2}},1)=\frac{1}{2\sqrt{3}} \sum _{c\ \in \mathcal{C}}
\left(\frac{c}{D_{T}} \right) _{3} \frac{1}{\wp \left(\frac{c\omega }{D} \right)
-1}+\frac{1}{3\sqrt{3}}\sum _{c\ \in \mathcal{C}}
\left(\frac{c}{D_{T}} \right) _{3}  . $$
This proves Theorem $1.6 $.
\par \vskip 0.4cm

  {\bf  Lemma  3.2.}
  $$ \sum \limits _{c\ \in \mathcal{C}}
\left(\frac{c}{D_{T}} \right) _{3}=\left\{
\begin{array}{l}  \sharp \mathcal{C}  \qquad  \hbox{if } \quad
T=\emptyset  ;\\
0  \qquad  \hbox{if } \quad T\neq \emptyset  .
\end{array}
\right.$$
\par \vskip 0.3cm

{\bf  Proof }$.\quad $  This could be verified by
the definition of cubic residue symbol( see [Ire-Ro] ).
   \par \vskip 0.4cm

{\bf  Lemma  3.3.}\quad
$ \sum \limits _{T} \left(\frac{c}{D_{T}} \right) _{3}=
\mu \cdot 2^{t},\quad $\\
 where $\ \mu \in \{ \pm 1,\pm \tau \pm \tau ^{2}\},\ $
$t=\sharp \left\{ \pi :\ \pi\equiv 1(mod\ 6)
 \hbox{ is a prime element,}\ \pi |D, \  \hbox{and }
\left({c}/{\pi } \right) _{3}=1 \right\},\ $
$c \in O_{K} $  and $D$ are relatively prime.
           \par \vskip 0.3cm

{\bf  Proof }$.\quad $By $\sum \limits _{T} \left(\frac{c}{D_{T}} \right) _{3}=
\left(1+\left(\frac{c}{\pi _{1}} \right) _{3} \right) \cdots
\left(1+\left(\frac{c}{\pi _{n}} \right) _{3} \right)
$ and the definition of cubic residue symbol, the lemma could
be verified easily.

 {\bf  Lemma 3.4.}
$$ (1)  \sum \limits _{T} 2^{n-t(T)} \left(\frac{c}
{D_{T}} \right) _{3}=(-\tau ^{2}) ^{t_{\tau }} \cdot 3^{t_{1}}
(1-\tau ) ^{t_{\tau }+t_{\tau ^{2}}} \quad ; $$
$$(2) \quad \sum \limits _{T} (-1) ^{t(T)} \left(\frac{c}
{D_{T}} \right) _{3}=\left\{
\begin{array}{l}  0  \qquad  \hbox{if }
\quad t_{\tau }+t_{\tau ^{2}}
< n \quad ; \\
(1-\tau ) ^{t_{\tau }+t_{\tau ^{2}}} \cdot (-\tau ^{2}) ^{t_{\tau ^{2}}}
 \qquad  \hbox{if } \quad t_{\tau }+t_{\tau ^{2}}=n \quad .
\end{array}
\right.$$
$$(3) \quad \sum \limits _{T}\sum \limits _{c \in \mathcal{C}} 2^{n-t(T)}
\left(\frac{c} {D_{T}} \right) _{3} =2^{n} \cdot
\sharp \mathcal{C}.$$
where $\ c \in O_{K},\ $  and $D$ are relatively prime,
the sum $\sum \limits _{T}$ is taken for $T$ runs over subsets
of  $\{1,\cdots , n\}$,
 $t=t(T)=\sharp T$( but $t=0$ when $T=\emptyset $),
 $$ t_{1}=\sharp \left\{\pi _{k}:\ \left(\frac{c}{\pi _{k}} \right) _{3}=1
\right\},\quad t_{\tau }=\sharp \left\{\pi _{k}:\ \left(\frac{c}{\pi _{k}} \right)
_{3}=\tau \right\},\quad t_{\tau ^{2}}=\sharp \left\{\pi _{k}:\ \left(\frac{c}
{\pi _{k}} \right) _{3}=\tau ^{2} \right\} $$
$t_{1}+t_{\tau }+t_{\tau ^{2}}=n \quad .$

{\bf  Proof }.\quad  Note that
$$\sum \limits _{T} 2^{n-t(T)} \left(\frac{c} {D_{T}} \right) _{3}=
\left(2+\left(\frac{c}{\pi _{1}} \right) _{3} \right) \cdots
\left(2+\left(\frac{c}{\pi _{n}} \right) _{3} \right) $$
$$\sum \limits _{T} (-1) ^{n-t(T)} \left(\frac{c} {D_{T}} \right) _{3}=
\left(1-\left(\frac{c}{\pi _{1}} \right) _{3} \right) \cdots
\left(1-\left(\frac{c}{\pi _{n}} \right) _{3} \right) .$$
Then by the definition of cubic residue symbol and Lemma
3.2, we could obtain the results.\par
\vskip 0.4cm

{\bf  Lemma  3.5 .} \quad  For the Weierstrass
$  \quad \wp -$function $\wp (z,L_{\omega })$ in Theorem
$4.3.1$   and any $c\ \in \mathcal{C},\ $  we have
$$ v_{3} \left(\wp \left(\frac{c\ \omega }{D},\ L_{\omega } \right) -1
\right)=\frac{1}{3}  . $$
\par \vskip 0.4cm

{\bf  Proof }.\quad  We need the following two lemmas
from [Ste ,\ P.128]:
          \par \vskip 0.3cm

 {\bf  Lemma  } $(C).\quad $Let $\Delta  $ be a square-free
integer in $K=\Q(\sqrt{-3})$ relatively prime to $\sqrt{-3}$,
$\beta \in O_{K}$ is relatively prime to  $\Delta  $ . And
$$ \varphi (\Delta  )=\left\{
\begin{array}{l}  \Delta  ^{2/(\Delta  \overline{\Delta  }-1)}
 \qquad \hbox{if }\ \Delta  \ \hbox{ is a prime };\\
1  \qquad  , \qquad  \hbox{ otherwise  }.
\end{array}
\right.$$
Then $\varphi (\Delta  )\wp (\beta \omega /\Delta  )$
is a unit.
        \par \vskip 0.3cm

 {\bf  Lemma  }$(D). \quad $ Let $r> 0$ be a rational integer.
Assume $\beta ,\ \Delta  ,$ and $  \gamma $ are integers of
$K=\Q(\sqrt{-3})$ ,  both $\Delta $ and $ \gamma $
are relatively prime to $\sqrt{-3}$, and$\beta $ is
relatively prime to $\Delta  $. Let $\lambda =
\frac{1}{2}(1-3^{1-r}),\ $ $\varphi (\Delta  )$ be as in
Lemma $(C)  .$ Then
$$ 3^{-\lambda }\varphi (\Delta  )\left\{\wp \left(\beta \omega /\Delta  \right)
-\wp \left(\gamma \omega /(\sqrt{-3}) ^{r} \right)
\right\}$$
is a unit.
          \par \vskip 0.3cm

Now we turn to the proof of Lemma  $ 3.5 . \quad $
In Lemma  $(D)$, let $r=2,\ \gamma =-1,\ $
$\Delta  =D,\ $ $\beta =c,\ $ then $\lambda
=\frac{1}{2}(1-3^{1-2})=\frac{1}{3},\ $ and we know
$$3^{-\frac{1}{3}} \varphi (D) \left\{\wp \left(
\frac{c\ \omega }{D} \right)-\wp \left(-\omega /
(\sqrt{-3}) ^{2} \right) \right\}=\theta $$
is a unit. $\varphi (D)=\varphi (\Delta  )$ is as in$(C)$.
Since  $\quad  \wp \left(-\omega /(\sqrt{-3}) ^{2} \right)=
\wp (\omega /3)=1,\ $ so
$$ \wp \left(\frac{c\ \omega }{D} \right)-1=\wp \left(\frac{c\ \omega }{D} \right)
-\wp \left(-\omega /(\sqrt{-3}) ^{2} \right)=3^{\frac{1}{3}} \varphi (D) ^{-1}
\theta $$
Thus we have  $v_{3}(\varphi (D))=0,\ $  and
 $$ v_{3} \left(\wp \left(\frac{c\ \omega }{D} \right)-1 \right)
=v_{3}(3^{\frac{1}{3}} \varphi (D) ^{-1} \theta )=\frac{1}{3}. $$
 This proves Lemma  3.5 .
          \par \vskip 0.4cm

{\bf  Proof of Theorem  1.7  .}\quad
For each subset $T$ of $\{1,\ \cdots ,\ n \}$,
multiply the two sides of formula $(1.6 )$ in Theorem $1.6 $
by $2^{n-t(T)},\ $  and then add them up, we obtain
\begin{align*}
&\sum _{T} 2^{n-t(T)}\frac{D}{\omega }\left(\frac{9}{D_{T}} \right) _{3}
L_{S}(\overline
{\psi } _{D_{T} ^{2}},1)\\
&=\frac{1}{2\sqrt{3}} \sum _{T} 2^{n-t(T)} \sum _{c\ \in \mathcal{C}}
\left(\frac{c}{D_{T}} \right) _{3} \frac{1}{\wp \left(\frac{c\ \omega }{D} \right)
-1}+\frac{1}{3\sqrt{3}} \sum _{T} 2^{n-t(T)}\sum _{c\ \in \mathcal{C}}
\left(\frac{c}{D_{T}} \right) _{3} \quad .
\end{align*}
Then by Lemma  $3.4(3),\ $  we have
\begin{align*}
&\sum _{T} 2^{n-t(T)}\frac{D}{\omega }\left(\frac{9}{D_{T}} \right) _{3}
L_{S}(\overline{\psi } _{D_{T} ^{2}},1)\\
&=\frac{1}{2\sqrt{3}} \sum _{c\ \in \mathcal{C}}\frac{1}
{\wp \left(\frac{c\ \omega }{D} \right) -1}\sum _{T} 2^{n-t(T)}\left(\frac{c}
{D_{T}} \right) _{3} +\frac{2^{n}}{3\sqrt{3}}\cdot \sharp \ \mathcal{C}
 \qquad  (3.18)
\end{align*}
By Lemma  3.4(1) we know,
$$ \quad \sum \limits _{T} 2^{n-t(T)} \left(\frac{c}
{D_{T}} \right) _{3}=(-\tau ^{2}) ^{t_{\tau }} \cdot 3^{t_{1}}
(1-\tau ) ^{t_{\tau }+t_{\tau ^{2}}}  \ ,  $$
Note that $1-\tau $ and $\sqrt{-3}$ are associated each other
(Both are prime elements in $O_{K}$), so
$v_{3}(1-\tau )=v_{3}(\sqrt{-3})=\frac{1}{2} \quad .$
Hence
$$ v_{3} \left(\quad \sum \limits _{T} 2^{n-t(T)} \left(\frac{c}
{D_{T}} \right) _{3} \right)=t_{1}+\frac{t_{\tau }+t_{\tau ^{2}}}{2}=
\frac{n+t_{1}}{2}\geq \frac{n}{2}$$
Also we have $\sharp \ \mathcal{C}=\prod \limits _{k=1} ^{n}(\pi _{k}
\overline{\pi _{k}}-1),\ $ so by the assumption we have
 $v_{3}(\sharp \ \mathcal{C})\geq n
,\ $  therefore
$$ v_{3} \left(\frac{2^{n}}{3\sqrt{3}}\cdot \sharp \ \mathcal{C} \right)
\geq n-\frac{3}{2}\geq \frac{n}{2}-1  \qquad (n\geq 1) $$
 so  by Lemma  $ 3.5 $ we know that the first term
 in the right side of $(3.18)$ has $3-$adic valuation
 $\geq\frac{n}{2}-\frac{1}{3}-\frac{1}{2}=
\frac{n}{2}-\frac{5}{6}> \frac{n}{2}-1\ .$
Therefore the right side of $(3.18)$ has $3-$adic valuation
$\geq \frac{n}{2}-1.\ $  Hence
$$ v_{3} \left(\sum _{T} 2^{n-t(T)}\frac{D}{\omega }
\left(\frac{9}{D_{T}} \right) _{3}
L_{S}(\overline{\psi } _{D_{T} ^{2}},1) \right)\geq
\frac{n}{2}-1  . $$
Also by Lemma  $ 3.1 $ we know that if
$\quad T=\{1,\ \cdots ,\ n \}$ then $L_{S}(\overline{\psi }
_{D_{T} ^{2}},1)=L(\overline{\psi } _{D ^{2}},1) \quad ;$
 and if $\quad T=\emptyset $ then
 $L_{S}(\overline{\psi }_{D_{T} ^{2}},1)
=L_{S}(\overline{\psi }_{1},1)=L(\overline{\psi } _{1},1)
\prod \limits _{k=1} ^{n} \left(1-\frac{1}{\pi _{k}} \right)$.
From  [Ste] we  know
 $$ L(\overline{\psi } _{1},1)=L(\psi _{1},\
1)=L_{1}(1)=(\frac{\sqrt{3}}{9})\omega  \ ,  $$
where $\psi _{1}$ is the Hecke characters of the  number field
$K=\Q(\sqrt{-3})$ corresponding to the elliptic curve
$\quad E:\ y^{2}=x^{3}-2^{4}3^{3}$ .
 Thus  $\quad L_{S}(\overline{\psi }_{1},1)=\frac{\sqrt{3}}{9}
 \omega \prod \limits _{k=1} ^{n} \left(1-\frac{1}{\pi _{k}}
 \right) .\ $ So we have
\begin{align*}
v_{3} \left(L_{S}(\overline{\psi } _{1},\ 1)/\omega \right)&=
v \left(\frac{\sqrt{3}}{9} \right)+\sum _{k=1} ^{n} v_{3}(\pi _{k}-1)\\
&\geq n-\frac{3}{2}\geq \frac{n}{2}-1 \quad (\hbox{ Since }v_{3}\ (\pi _{k}-1)
\geq 1)
\end{align*}
Now we use induction method on $n$ to prove
$\quad v_{3} \left(L (\overline{\psi } _{D^{2}} ,\ 1)/
 \omega \right)\geq \frac{n}{2}-1 $ .
when $n=1$, $D=\pi _{1},\ $  and$L_{S}(\overline{\psi } _{1},
\ 1)=\frac{\sqrt{3}}{9}\omega \left(1-\frac{1}{\pi _{1}}
\right).\ $
Since $\quad v_{3}(\pi _{1}-1)\geq 1,\ $ so
$$v_{3}\left(L_{S}(\overline{\psi } _{1},\ 1) /\omega \right)
\geq 1-\frac{3}{2}=-\frac{1}{2}  . $$
Also we have
$$ v_{3}\left(2\frac{\pi _{1}}{\omega }\left(\frac{9}{D_{\emptyset }}
\right) _{3}L_{S}(\overline{\psi } _{1},\ 1)+\frac{\pi _{1}}{\omega }
\left(\frac{9}{\pi _{1}}
\right) _{3}L(\overline{\psi } _{\pi _{1} ^{2}},\ 1) \right)\geq
\frac{1}{2}-1=-\frac{1}{2} $$
  Therefore
$$ v_{3}\left(L(\overline{\psi } _{\pi _{1} ^{2}},\ 1)/\omega  \right)=
v_{3} \left(\frac{\pi _{1}}{\omega }\left(\frac{9}{\pi _{1}}
\right) _{3}L(\overline{\psi } _{\pi _{1} ^{2}},\ 1) \right)\geq -\frac{1}{2}
=\frac{1}{2}-1  . $$
Assume our conclusion is true for $1,\ 2,\ \cdots ,\ n-1$,
and consider the case $n$, $D=\pi _{1} \cdots \pi _{n}  .$
For any nonempty subset $T$ of $\{1,\ \cdots , n\}$, put
$t=t(T)=\sharp \ T$, by  Lemma  $ 3.1 $ we know
$$ 2^{n-t(T)}\ \frac{D}{\omega }\left(\frac{9}{D_{T}} \right) _{3}
L_{S}(\overline{\psi } _{D_{T} ^{2}},1) =  2^{n-t(T)}\ \frac{D}{\omega }
\left(\frac{9}{D_{T}} \right) _{3}L(\overline{\psi } _{D_{T} ^{2}},1)
\prod _{\pi _{k}|\widehat{D} _{T} } \left(1-\left(\frac{D_{T}}
{\pi _{k}} \right) _{3}\frac{1}{\pi _{k}} \right) $$
 Hence
\begin{align*}
v_{3} \left(2^{n-t(T)}\ \frac{D}{\omega }\left(\frac{9}{D_{T}} \right) _{3}
L_{S}(\overline{\psi } _{D_{T} ^{2}},1) \right)&=v_{3} \left(L_{S}
(\overline{\psi } _{D_{T} ^{2}},1)/\omega  \right)\\
&=v_{3} \left(L
(\overline{\psi } _{D_{T} ^{2}},1)/\omega  \right)
+\sum _{\pi _{k}|\widehat{D} _{T}} v_{3} \left(1-\left(\frac{D_{T}}
{\pi _{k}} \right) _{3}\frac{1}{\pi _{k}} \right)
\end{align*}
Note that $T$ is a proper subset, by induction assumption
 we have  $v_{3} \left(L(\overline{\psi }
_{D_{T} ^{2}},1)/\omega  \right)\geq \frac{t(T)}{2}-1  .$
Since $\quad \left(\frac{D_{T}}{\pi _{k}} \right) _{3}=1,\
  \tau ,\ \hbox{or } \ \tau ^{2}$
  (since $\pi_{k}|\widehat{D} _{T}$), so
\begin{align*}
v_{3} \left(1-\left(\frac{D_{T}}
{\pi _{k}} \right) _{3}\frac{1}{\pi _{k}} \right)=v_{3}(\pi _{k}-1)\ ,\
v_{3}(\pi _{k}-\tau )\ ,\ \hbox{or } \  v_{3}(\pi _{k}-\tau ^{2})\geq
1\ ,\ \frac{1}{2}\ ,\ \hbox{or }\ \frac{1}{2} \quad .
\end{align*}
thus
$$ \sum _{\pi _{k}|\widehat{D} _{T}} v_{3} \left(1-
\left(\frac{D_{T}}
{\pi _{k}} \right) _{3}\frac{1}{\pi _{k}} \right)\geq
\frac{n-t(T)}{2}. $$
  Therefore
$$ v_{3} \left(2^{n-t(T)}\ \frac{D}{\omega }\left(\frac{9}{D_{T}} \right) _{3}
L_{S}(\overline{\psi } _{D_{T} ^{2}},1) \right)\geq \frac{t(T)}{2}-1+
\frac{n-t(T)}{2}=\frac{n}{2}-1  . $$
Also  when   $ T=\emptyset $, we above have proved
 $$ v_{3} \left(2^{n}\ \frac{D}{\omega }\left(\frac{9}{D_{\emptyset }} \right)
_{3} L_{S}(\overline{\psi } _{1},1) \right)=v_{3} \left(L_{S}
(\overline{\psi } _{1} ,\ 1)/\omega \right)\geq \frac{n}{2}-1  . $$
  Therefore
\begin{align*}
v_{3} \left(L(\overline{\psi } _{D ^{2}},\ 1)  / \omega  \right)\\
&=
v_{3} \left(2^{n-n}\ \frac{D}{\omega }\left(\frac{9}{D} \right) _{3}
L_{S}(\overline{\psi } _{D ^{2}},1) \right)\\
&=v_{3} \left( \left(
\sum _{T} 2^{n-t(T)}\ \frac{D}{\omega }\left(\frac{9}{D_{T}} \right) _{3}
L_{S}(\overline{\psi } _{D_{T} ^{2}},1) \right)\right.\\
&\left.\mbox{}\hspace{2cm}-\left(
\sum _{T\subsetneqq \{ 1,\cdots , n\}} 2^{n-t(T)}\ \frac{D}{\omega }
\left(\frac{9}{D_{T}} \right) _{3}
L_{S}(\overline{\psi } _{D_{T} ^{2}},1) \right) \right)\\
&\geq \frac{n}{2}-1
\end{align*}
This proves our conclusion by induction, and completes the
proof of Theorem $ 1.7  $ .
\vskip 0.4cm

{\bf  Proof of Theorem 1.8 }$.\quad $
Since $\quad \pi _{k}\equiv 1 (mod\ 6),\ $ $\sharp \ \mathcal{C}=
\prod \limits _{k=1} ^{n} (\pi _{k} \overline{\pi }_{k}-1)\equiv 0 (mod\
6),\ $ so we may choose the set $\mathcal{C},\ $
in such a way that  $\pm c,\ \pm \tau c,\ \pm \tau ^{2}c
\in \mathcal{C}$( when $c \in \mathcal{C}$). That is,
 when  $c \in \mathcal{C}$, then all its associated
 elements are in $\mathcal{C} \quad .$
 Let $\quad  V=\{c \in \mathcal{C}:\ c\equiv 1(mod\ 3)\},\ $
 Then $$\quad \mathcal{C}=\bigcup \limits _{\mu \
 \in \ \{\pm 1, \pm \tau , \pm \tau ^{2}\}} \mu V \quad .$$
Since
$$\quad \left(\frac{-c}{D_{T}} \right) _{3}=\left(\frac{c}
{D_{T}} \right) _{3}  \ ,  \quad  \left(\frac{\tau }{D_{T}} \right) _{3}
=1,\ \tau ,\ \hbox{or } \ \tau ^{2} \quad , $$
$$ \wp \left( \frac{-c\ \omega }{D} ,\ L_{\omega } \right)=
\wp \left( \frac{c\ \omega }{D} ,\ L_{\omega } \right) \quad ,$$
$$ \wp (\tau z,\ L_{\omega })=\wp (\tau z,\ \tau L_{\omega })=
\frac{1}{\tau ^{2}} \wp (z,\ L_{\omega })=\tau \wp (z,\ L_{\omega })
\quad , $$
$$ \wp (\tau ^{2} z,\ L_{\omega })=\wp (\tau ^{2} z,\ \tau ^{2}
L_{\omega })=\frac{1}{\tau ^{4}} \wp (z,\ L_{\omega })
=\tau ^{2} \wp (z,\ L_{\omega }) \quad , $$
therefore
\begin{align*}
S^{*}(D)&=\frac{1}{\sqrt{3}} \ \sum _{c\ \in V}
\frac{1}{\wp (\frac{c\ \omega }{D},\ L_{\omega })-1}
\sum _{T} 2^{n-t(T)}\ \left( \frac{c}{D_{T}} \right) _{3}\\
&+ \frac{1}{\sqrt{3}} \ \sum _{c\ \in V}
\frac{1}{\wp (\frac{\tau \ c\ \omega }{D},\ L_{\omega })-1}
\sum _{T} 2^{n-t(T)}\ \left( \frac{\tau \ c}{D_{T}} \right) _{3}\\
&+ \frac{1}{\sqrt{3}} \ \sum _{c\ \in V}
\frac{1}{\wp (\frac{\tau ^{2} \ c\ \omega }{D},\ L_{\omega })-1}
\sum _{T} 2^{n-t(T)}\ \left( \frac{\tau ^{2} \ c}{D_{T}} \right) _{3}\\
& =\frac{1}{\sqrt{3}} \sum _{T} 2^{n-t(T)} \left\{ \sum _{c\ \in V}
W_c\left(\frac{c}{D_{T}} \right) _{3} \right\}
\end{align*}
where
$$W_c=\frac{1}{\wp (\frac{c\ \omega }{D},\ L_{\omega })-1}+
\frac{\left(\frac{\tau }{D_{T}} \right) _{3}}
{\tau \ \wp (\frac{c\ \omega }{D},\ L_{\omega })-1}+
\frac{\left(\frac{\tau }{D_{T}} \right) _{3} ^{2}}
{\tau ^{2} \ \wp (\frac{c\ \omega }{D},\ L_{\omega })-1}.$$
Denote
$$ U(\wp )=\left( \wp \left( \frac{c\ \omega }{D} \right)-1 \right)
\left(\tau \ \wp \left( \frac{c\ \omega }{D} \right)-1 \right)
\left(\tau ^{2}\ \wp \left( \frac{c\ \omega }{D} \right)-1 \right) $$
then
$$ \frac{1}{\wp (\frac{c\ \omega }{D},\ L_{\omega })-1}+
\frac{\left(\frac{\tau }{D_{T}} \right) _{3}}
{\tau \ \wp (\frac{c\ \omega }{D},\ L_{\omega })-1}+
\frac{\left(\frac{\tau }{D_{T}} \right) _{3} ^{2}}
{\tau ^{2} \ \wp (\frac{c\ \omega }{D},\ L_{\omega })-1}$$
$$  = \frac{1}{U(\wp )} \cdot
\left\{
\begin{array}{c}  3  \qquad  \hbox{if } \quad \left(\frac{\tau }
{D_{T}} \right) _{3}=1 \quad ;\\
3\wp ^{2}  \qquad  \hbox{if } \quad \left(\frac{\tau }
{D_{T}} \right) _{3}=\tau  \quad ;\\
3\wp   \qquad  \hbox{if } \quad \left(\frac{\tau }
{D_{T}} \right) _{3}=\tau ^{2} \quad .
\end{array}
\right. $$
Also put
$$ V(\wp )=\left\{
\begin{array}{c}  1  \qquad  \hbox{if } \quad \left(\frac{\tau }
{D_{T}} \right) _{3}=1 \quad ;\\
\wp ^{2}  \qquad  \hbox{if } \quad \left(\frac{\tau }
{D_{T}} \right) _{3}=\tau  \quad ;\\
\wp   \qquad  \hbox{if } \quad \left(\frac{\tau }
{D_{T}} \right) _{3}=\tau ^{2} \quad .
\end{array}
\right. $$
where $\quad \wp =\wp \left( \frac{c\ \omega }{D} \right), $
then
\begin{align*}
S^{*}(D)&=\frac{1}{\sqrt{3}} \sum _{T} 2^{n-t(T)} \left\{ \sum _{c\ \in V}
\frac{3V(\wp )}{U(\wp )} \cdot \left(\frac{c}{D_{T}} \right) _{3}
\right\}\\
&=\sqrt{3} \sum _{c\ \in V} \frac{V(\wp )}{U(\wp )} \sum _{T} 2^{n-t(T)}
\left(\frac{c}{D_{T}} \right) _{3}
\end{align*}
Since
$$ \tau \ \wp \left( \frac{c\ \omega }{D} \right)-1=
\left( \wp \left( \frac{c\ \omega }{D} \right)-\tau ^{2} \right)/\tau
^{2}=\frac{1}{\tau ^{2}} \left(\left( \wp \left( \frac{c\ \omega }{D} \right)
-1 \right)+(1-\tau ^{2}) \right)  \ ,  $$
$$ \tau ^{2} \ \wp \left( \frac{c\ \omega }{D} \right)-1=
\frac{1}{\tau } \left( \wp \left( \frac{c\ \omega }{D} \right)-\tau  \right)
=\frac{1}{\tau }\left( \left( \wp \left( \frac{c\ \omega }{D} \right)-1 \right)
+(1-\tau ) \right)  \ ,  $$
so by Lemma  $ 3.5 $ we obtain
$$ v_{3} \left( \tau \ \wp \left( \frac{c\ \omega }{D} \right)-1 \right)=
v_{3} \left( \tau ^{2}\ \wp \left( \frac{c\ \omega }{D} \right)-1 \right)
=v_{3} \left( \wp \left( \frac{c\ \omega }{D} \right)-1 \right)
=\frac{1}{3}  \ ,  $$
thus
$$ v_{3}(U(\wp ))=3v_{3} \left( \wp \left( \frac{c\ \omega }{D} \right)-1
\right)=3\cdot \frac{1}{3}=1 \quad . $$
Also we have
$$ v_{3} \left( \sum _{T} 2^{n-t(T)} \left( \frac{c}{D_{T}} \right) _{3}
\right) \geq \frac{n}{2} \quad (\hbox{ see Lemma  }3.4(1)) $$
  therefore
$$ v_{3}(S^{*}(D))\geq \frac{1}{2}+v_{3} \left( \frac{V(\wp )}{U(\wp )}
\right)+v_{3} \left( \sum _{T} 2^{n-t(T)} \left( \frac{c}{D_{T}} \right)
_{3} \right)\geq \frac{1}{2}-1+\frac{n}{2}=\frac{n-1}{2} \quad . $$
 This proves the Proposition .

\vskip 1cm
\begin{center}{\bf {\large \bf R}eferences}\end{center}
 \baselineskip 0pt
\parskip 0pt
\begin{description}

\item[[Bir-Ste]] B. J. Birch and N. M.Stephens, The parity of
the rank of the Mordell-Weil
group, Topology 5 (1996), 295-299.

\item[[B-SD]]   B. J. Birch and H.P.F.Swinnerton-Dyer, Notes
on elliptic curves II,
J. Reine Angew . Math. 218(1965), 79-108.

\item[[Co-Wi]] J. Coates, A. Wiles, On the conjecture of Birch
and  Swinnerton-Dyer, Invent.
Math., 39(1977), 223-251.

\item[[Go-Sch]] C.Coldstein and N.Schappacher,
S$\acute{e}$ries d' Eisenstein et fonction
L de courbes ellipliques $\grave{a}$ multiplication complexe,
 J. Reine Agew . Math.,
327(1981), 184-218.

\item[[Ire-Ro]]  K.Ireland  and M.Rosen,   A Classical
Introduction to Modern Number Theory, GTM 84, Springer-Verlag,
New York, 1990.

\item[[Law]] D.F.Lawden,  Elliptic Functions and Applications,
Applied Mathematical Sciences Vol.80, Springer-Verlag,
New York, 1989.

\item[[Razar]]

\item[[Ru\ 1]] K.Rubin, Tate-Shafarevich groups and
L-functions of elliptic curves with complex
multiplication, Invent. Math., 89(1987), 527-560.

\item[[Ru\ 2]] K.Rubin , The `` main conjectures ''
of Iwasawa theory for imaginary quadratic
fields, Invent.  Math., 103(1991), 25-68.

\item[[Sil\ 1]] J. H. Silverman,  ``  The Arithmetic of
Elliptic Curves '', GTM 106,
Springer-Verlag, New York, 1986.

\item[[Sil\ 2]] J. H. Silverman, ``  Advanced Topics in the
Arithmetic of  Elliptic Curves '',  GTM 151, Springer-Verlag,
1994.

\item[[Ste]] N. M. Stephens, The diophantine equation
$x^{3}+y^{3}=Dz^{3}$ and the conjectures
of Birch and Swinnerton-Dyer, J. Reine Angew. Math.,
231(1968), 121-162

\item[[Tun]] J. B. Tunnell , A classical Diophantine problem
and modular forms of weight
$\frac{3}{2}$. Invent. Math., 72(1983), 323-334.

\item[[Wa]] L. C. Washington, `` Introduction to Cyclotomic
Fields '', Springer-Verlag, New York, 1982

\item[[We]] A. Weil, Elliptic functions according to
Eisenstein and Kronecker, Springer, 1976

\item[[Zhang]] ZHANG Xianke, Introduction to Algebraic Number
Theory,  Hunan Edu. Press, Hunan, China, 1999.
\item[[Zhao]]   ZHAO Chunlai, A criterion for elliptic curves
with lowest 2-power in L(1),
Math. Proc. Cambridge Philos. Soc. 121(1997), 385-400.
  \end{description}
  \par \vskip 0.6cm

Tsinghua University, The Center for Advanced Study,
\par \vskip 0.06cm
Beijing 100084, R. R. China
 \par \vskip 0.4cm

Tsinghua University, Department of Mathematical Sciences,
\par \vskip 0.06cm
Beijing 100084, P. R. China
\par \vskip 0.14cm

xianke@tsinghua.edu.cn

(This paper was published in :

Acta Arithmetica, 103.1(2002), 79-95 );

Manuscripta Math. 108, 385-397 (2002)

\end{document}